\documentclass[article,11pt]{article}
\usepackage{graphicx}
\usepackage{amsmath,amsthm}
\usepackage{mathrsfs,amsfonts}
\usepackage{amsfonts}
\usepackage{amsmath,amssymb,amsfonts}
\usepackage{amssymb}
\usepackage{bm}
\usepackage{color}
\usepackage{amsmath}  
\usepackage{amsfonts} 
\usepackage{graphicx} 
\usepackage{braket}
\usepackage{enumerate}
\usepackage{indentfirst}
\usepackage{longtable}
\usepackage{setspace}
\allowdisplaybreaks[4]

\def\al{{\alpha}}\def\be{{\beta}}
\def\ep{{\epsilon}}\def\ga{{\gamma}}
\def\la{{\lambda}}\def\si{{\sigma}}

\def\<{\left<}\def\>{\right>}\def\({\left(}\def\){\right)}

\newfam\msbmfam\font\tenmsbm=msbm10\textfont
\msbmfam=\tenmsbm\font\sevenmsbm=msbm7

\scriptfont\msbmfam=\sevenmsbm\def\bb#1{{\fam\msbmfam #1}}

\def\EE{\bb E}\def\FF{\bb F}\def\HH{\bb H}

\def\RR{\bb R}\def\SS{\bb S}
\def\cB{{\cal B}}\def\cL{{\cal L}}
\def\cF{{\cal F}}
\def\cH{{\cal H}}
\def\cL{{\cal L}}
\def\cP{{\cal P}}\def\cS{{\cal S}}\def\cT{{\cal T}}\def\cU{{\cal U}}\def\cQ{{\cal Q}}  \def\cQ{{\cal Q}}

\usepackage{amsmath,amsfonts,amsthm,amssymb}    
\usepackage{gensymb, mathrsfs}                  

\usepackage{accents}
\DeclareMathSymbol{\widehatsym}{\mathord}{largesymbols}{"62}

\def\*#1{\mathbf{#1}}

\newcommand{\tr}{\operatorname{tr}}

\newcommand{\Ker}{\operatorname{Ker}}

\renewcommand{\bar}{\overline}
\renewcommand{\tilde}{\widetilde}
\renewcommand{\phi}{\varphi}

\renewcommand{\qed}{\hfill \ensuremath{\Box}}

\makeatletter \renewcommand\d[1]{\ensuremath{%
  \;\mathrm{d}#1\@ifnextchar\d{\!}{}}}
\makeatother

\theoremstyle{plain}
\newtheorem{thm}{Theorem}[section]
\newtheorem{lem}[thm]{Lemma}

\theoremstyle{definition}


\newcommand{\beq}{\begin{equation}}
\newcommand{\eeq}{\end{equation}}

\definecolor{c}{rgb}{0.9,0.3,0.1}
\definecolor{b}{rgb}{0.1,0.3,0.9}

\newtheorem{remark}{Remark}[section]
\newtheorem{lemma}{Lemma}[section]
\newtheorem{theorem}{Theorem}[section]

\renewcommand{\theequation}{\arabic{section}.\arabic{equation}}

\def\al{{\alpha}}\def\be{{\beta}}
\def\ep{{\epsilon}}
\def\ga{{\gamma}}\def\la{{\lambda}}
\def\si{{\sigma}}

\def\<{\left<}\def\>{\right>}\def\({\left(}\def\){\right)}

\newfam\msbmfam\font\tenmsbm=msbm10\textfont
\msbmfam=\tenmsbm\font\sevenmsbm=msbm7
\scriptfont\msbmfam=\sevenmsbm\def\bb#1{{\fam\msbmfam #1}}

\def\EE{\bb E}\def\HH{\bb H}
\def\RR{\bb R}

\def\cB{{\cal B}}\def\cF{{\cal F}}
\def\cH{{\cal H}}\def\cK{{\cal K}}\def\cL{{\cal L}}\def\cP{{\cal P}} 
\def \cS{{\cal S}}\def \cR{{\cal R}}\def \cN{{\cal N}}\def \cO{{\cal O}}

\numberwithin{equation}{section}

\begin{document}

\baselineskip=18pt

\title{\large \bf  General  mean-field stochastic linear quadratic control problem   driven by L\'evy processes with random coefficients}
\author{Yanyan Tang\thanks{Department of Mathematics, Southern University of Science and Technology, Shenzhen 518055, China ({\tt 12131233@mail.sustech.edu.cn}). } \;
and
Jie Xiong\thanks{ Department of Mathematics and  Shenzhen International Center for Mathematics, Southern University of Science and Technology, Shenzhen 518055, China ({\tt xiongj@sustech.edu.cn}).
}
}
\date{}
\maketitle
 \bigskip

\noindent \textbf{Abstract.} 
This paper studies a stochastic mean-field linear-quadratic optimal control problem with random coefficients. The state equation is a general
linear stochastic differential equation with mean-field terms $\EE X(t)$ and $\EE u(t)$ of the state and the control processes and is
 driven by a Brownian motion and a Poisson random measure. By the   coupled system of Riccati equations,  an explicit expressions for the optimal state feedback control is obtained.  As a by-product, the non-homogeneous stochastic linear-quadratic control problem with random coefficients and L\'evy driving noises is also studied.
\bigskip


\section{Introduction} \label{sec1}
\setcounter{equation}{0}
\renewcommand{\theequation}{\thesection.\arabic{equation}}
 Let $(\Omega, \mathcal{F}, \mathbb{P})$ be a probability space with a filtration $\FF\equiv(\mathcal{F}_{t})_{t\geq0}$ satisfying the usual conditions. In this space, there is an one-dimensional standard Brownian motion  $W$ and an independent Poisson random measure $N(dt,dz)$ with intensity  $\nu(dz)$ which is a $\sigma$-finite measure on a measurable space $(U_0,\cU_0)$. We write $\tilde{N}(dt,dz)=N(dt,dz)-\nu(dz) dt$ for the compensated  martingale random measure. 

 We consider the  following linear stochastic  mean-field system driven by the  Brownian motion  $\{W(t)\}$ and the Poisson random
 measure $\{{N}(dt,dz)\}$:
\begin{equation}\label{sate}
\left\{\begin{array}{ccl}
dX(t)&=&\big(A(t)X(t)+A_1(t)\EE[X(t)]+B(t)u(t)+B_1(t)\EE[u(t)]\big)dt\\
&&\ \ +\big(C(t)X(t)+C_1(t)\EE[X(t)]+D(t)u(t)+D_1(t)\EE[u(t)] \big)dW(t)\\
&&\ \ +\int_{U_0}\big(\al(t,z)X(t-)+\al_1(t,z)\EE[X(t-)]+\be(t,z)u(t)+\be_1(t)\EE[u(t)]\big)\tilde{N}(dt,dz)\\
X(0)&=&x, \ t\in[0,T],
\end{array}\right.
\end{equation}
where the initial state $x\in\RR^n$, and the coefficients  $A, A_1, C, C_1:[0,T]\times\Omega\to\RR^{n\times n}$, $B, B_1, D, D_1:[0,T]\times\Omega\to\RR^{n\times m}$, $\al(\cdot,\cdot), \al_1(\cdot,\cdot):[0,T]\times U_0 \times\Omega \to\RR^{n\times n}$ and $\be(\cdot,\cdot), \be_1(\cdot,\cdot):[0,T]\times U_0\times\Omega\to , \RR^{n\times m}$ are measurable mappings which are $\cF_t$-adapted,   $u$ is the control process. Note that the process $X=X^u$ depends on the control $u$ and we omit the superscript $u$ for convenience.

An admissible control $u$  is defined as an   $\{\mathcal{F}_t, \ 0 \leq t\leq T\}$ adapted process taking values in $U\subset\RR^m$  such that 
\begin{eqnarray} \label{uu}
\mathbb{E}\left(\int_0^T |u(t)|^2 dt\right) < \infty.
\end{eqnarray} 
The set of all admissible controls  is denoted by $\cU_{ad}$.  Note that $\cU_{ad}$ is a Hilbert space. 
 For any admissible control $u\in\cU_{ad}$, we consider the following quadratic cost functional 
\begin{eqnarray}\label{cost1}
J(u) &=&\EE\bigg\{\big<GX(T), X(T)\big>
+\int_0^T \Big(\big<Q(t)X(t), X(t)\big> +  \big<Q_1(t)\EE[X(t)], \EE[X(t)]\big>\nonumber\\
&&\ \ \qquad \qquad\qquad +\big<R(t)u(t), u(t)\big>+\big<R_1(t)\EE[u(t)], \EE[u(t)]\big>\Big)dt\bigg\}
\end{eqnarray}
where $Q, Q_1:[0,T]\times\Omega\to\RR^{n\times n}$, $R, R_1:[0,T]\times\Omega\to\RR^{m\times m}$ are $\cF_t$-adapted processes and $G$ is an  $\cF_T$-measurable random matrix.

The objective of  the mean-field stochastic linear-quadratic (MF-SLQ in short) optimal control problem  with random coefficients and jumps is to choose an admissible control such that the cost functional (\ref{cost1}) is minimized. 
\\
{\bf{Problem (MF-SLQ)}}: 
 Find a control $u^*\in\cU_{ad}$  such that 
 \[J(u^*)=\inf_{u\in\cU_{ad}} J(u).\]
 If such a process $u^*$ exists, it is called an optimal control 
    and the corresponding state process $X^*$ is called the  optimal state process. We also refer to $(X^*, u^*)$ as the optimal pair.
\medskip

Mean-field stochastic differential equations (MF-SDEs) were originally used to model physical systems with a large number of interacting particles.
 In the dynamics of MF-SDEs, the complexity is reduced by substituting the interactions among all particles with their average. 
 Over the past decade, since Buchdahn et al \cite{{BR20091}, {BR20092}} and Carmona-Delarue \cite{{CR20131}, {CR20132}} 
 introduced  mean-field backward stochastic differential equations (MF-BSDEs) and mean-field forward-backward stochastic differential 
 equations (MF-FBSDEs), mean-field stochastic control problems and stochastic differential game problems have become a popular topic.  
 For some recent research see \cite{{BRJ2011},  {DB2015}, {GX}, {GXZ}, {HWX}, {Li2012},  {SWX}, {TX}, {WZ2022}, {ZSX}}. 
 MF-SLQ optimal control problem has attracted significant research interests. Firstly, it provides an elegant approach to solving the optimal control problem for continuous-time systems. By solving the Riccati equation, SLQ optimal control can derive the optimal state feedback controller, enabling the system to achieve optimal performance under given performance metrics. The  MF-SLQ problem was first studied by Yong \cite{Y2013}, who gave a specific expression for the optimal control $u$ with feedback form of the state process $X$, for the case where all coefficients are deterministic. After that, other  linear-quadratic mean field control and game problems have 
 been studied by  \cite{{Lsy2021}, {LWX}, {LSX}, {SY2020}, {SJ2017}, {WG2015}}.

In fact, as early as the 1960s, the application of SLQ optimal control problems with Poisson jumps has received some development (see \cite{ {JC1992}, 
 {KRA1961}, {sworder1969feedback}}). Recently, there has been a sharp increase in interest in  mean-field optimal control problems with jumps.  
 Shen and Siu \cite{ST2013} established a necessary and sufficient stochastic maximum principle for a mean-field model with randomness described by Brownian motion and Poisson jumps; Yang et al. \cite{YM2014}  generalised the model of \cite{ST2013} and discussed a stochastic optimal control problem with delay of mean field types, two sufficient maximum principles and one necessary maximum principle are established, and the mean-field type of game problem with jumps has also been studied by Moon \cite{M2022}.
  
 After MF-SLQ control problems with deterministic coefficients are investigated,  a natural question is what would happen if the coefficients are random. Since the Riccati equation corresponding to the SLQ problem with random coefficients (SRE in short form) is a highly nonlinear BSDE, proving the uniqueness of its solution is an extremely challenging task. In fact,  the special case where the SRE generator depends only linearly on the unknowns 
  was first studied by Bismut \cite{BI1976} by constructing a contraction mapping to make use of the fixed point theorem. For the general case,  Tang \cite{T2003} was able to solve it using a stochastic flow approach, i.e., the existence and uniqueness of the SRE solution was proved in the general case corresponding to a linear quadratic problem with random coefficients and state and control-dependent noise.  The LQ problem with jumps and  random coefficients has been studied by Meng \cite{Meng2014} and Zhang  et al. \cite{ZDM2020}. In \cite{Meng2014} Meng established the connection between the multidimensional  stochastic Riccati equation with jumps (SREJ for short) and the SLQ  problem and the related Hamiltonian system, and demonstrated the existence and uniqueness results of the multidimensional SREJ. We emphasize that the works \cite{Meng2014} and \cite{ZDM2020}  play an important role in our current study.

MF-SLQ control problem with random coefficients was posed as an open problem by Yong \cite{Y2013}. 
The main difficulty caused by random coefficients is the decoupling of the terms of the form $\EE[A(t)^\top Y(t)]$ according to $\EE[Y(t)]$, where $A(t)$ is a random coefficient 
and $Y(t)$ is a suitable adjoint process. To overcome this hurdle, Xiong and Xu \cite{XW2024} invented an extended LaGrange multiplier (ELM in short form)
 method 
to deal with MF-SLQ with random coefficients, where the state equation is driven by a Brownian motion.

Note that the state equation of  \cite{XW2024} does not involve the $\EE[u(t)]$ term. If we apply the method of \cite{XW2024} directly,
we would obtain an equation of the form $Ru+R_1\EE[u]+Q=0$, which is very hard to solve because of the 
randomness of $R$. In this article, we extend the method of  \cite{XW2024} to the case of MF-SLQ with random coefficients for state equation 
driven not only by a Brownian motion but also by a Poisson random measure. This is one of the main contributions of this article. 
The second contribution is the inclusion of the terms involving $\EE[u(t)]$ which makes the linear state equation in its most general form. Thirdly, this article 
makes contribution to the SLQ control problem with random coefficients in the following sense. The problem without the jumping Poisson random measure was
studied by Tang \cite{T2003} for the homogenous case; while that with  jumping Poisson random measure was studied by Meng \cite{Meng2014}.
Recently, Hu et al \cite{HSX} studied the non-homogeneous case driven by Brownian motion using BMO martingale techniques. In this paper, we study the non-homogeneous problem with jumps. More specifically, we establish the existence and uniqueness of the solution to the linear
BSDE with possibly unbounded coefficients without using the BMO theory.

The main purpose of this paper is to discuss in detail the MF-SLQ control problem with random coefficients, where the linear system is a multidimensional stochastic differential equation driven by  a one-dimensional Brownian motion and a Poisson random measure.
We first establish the stochastic maximum principle for 
Problem (MF-SLQ), including the existence and uniqueness of the optimal control. More importantly, 
we avoid the  form of $Ru+R_1\EE[u]+Q=0$, which is very difficult to solve  for the optimal control $u$,
 by making $\EE[u]=b$ a constraint and by using the ELM method.  Finally, we show that the optimal control has the state feedback representation.

The rest of the paper is organised as follows. In Section \ref{sec2},  we introduce 
 some useful notation and some existing results on MF-SDEs  and  MF-BSDEs driven by Poisson random  measure.  In Section \ref{sec3}, we present the main results of this paper.   In Section \ref{sec4},  we prove that the  Problem (MF-SLQ) has a unique optimal control, and  it satisfies an optimality system.   In Section \ref{sec5},  we solve the  control problem under the constraints $\EE [X]=a $ and $\EE [u]=b $ using the ELM method.    The SLQ control problem with operator coefficients and deterministic control of the frozen mean-field is studied in Section \ref{sec6}.  Throughout this article, $K$  is a positive  constant whose value can  be different from place to place.

\section{ Preliminaries } \label{sec2}
\setcounter{equation}{0}
\renewcommand{\theequation}{\thesection.\arabic{equation}}

In this section, we   present two preliminary theorems on the existence and uniqueness of solutions to  SDEs and BSDEs with jumps.   For the convenience of later uses, we introduce some notations. Let $\bar{X}$ be the expectation of the random variable $X$.  For any Euclidean space  $\HH=\RR^n, \RR^m, \SS^n_{+}$ (with $\SS^n_{+}$ being the set of all $n\times n$ positive semi-definite matrices)
 we introduce the following spaces:
\begin{itemize}  
\item $L^2(0,T; \HH)$: the space of all  deterministic  $\HH$-valued  functions with $\int_0^T\|x(t)\|_\HH^2dt<\infty.$
\item $L^2_{\cF_t}(\Omega; \HH)$: the space of all   $\HH$-valued $\cF_t$-measurable random variables $\xi$ with $\EE\|\xi\|_\HH^2<\infty.$ 
    \item $S^p_\FF(0,T; \HH)$: the space of all $\HH$-valued   $\cF_t$-adapted c\`{a}dl\`{a}g processes $X:[0,T]\times \Omega\to \HH$ satisfying 
   $$\EE\sup_{0\leq t\leq T} \|X(t,\omega)\|_\HH^p<\infty.$$
  \item $L_\FF^{2,p}(0,T; \HH)$: the space of all  $\HH$-valued   $\cF_t$-adapted c\`{a}dl\`{a}g processes  $X:[0,T]\times \Omega\to \HH$   satisfying
$$\EE\Big(\int_0^T\|X(t,\omega)\|_{\HH}^2dt\Big)^{p}<\infty.$$
   Especially, we denote $L_{\FF}^2(0,T; \HH)\equiv L_{\FF}^{2,1}(0,T; \HH)$.
   \item $L_\FF^\infty(0,T; \HH)$:  the space of all  $\HH$-valued   $\cF_t$-adapted c\`{a}dl\`{a}g  bounded processes.
   \item $L_{\cF_t}^\infty( \HH)$:  the space of all  $\HH$-valued   $\cF_t$-measurable  bounded  random variables.
    \item $\cL^2(U_0; \HH)$: the space of $\HH$-valued  measurable functions  $r=\{r(z),z\in U_0\}$ defined on the measure space $(U_0, \cB(U_0); \nu)$ satisfying
   $$\int_{U_0}\|r(z)\|_{\HH}^2\nu(dz)<\infty.$$

  \item $\cL_\FF^{2, p}([0,T]\times U_0; \HH)$: the space of $\cL^2(U_0; \HH)$-valued  and $\cF_t$-predictable processes  $r=\{r(t,\omega,z),z\in U_0\}$  satisfying
$$\EE\Big(\int_0^T\int_{U_0}\|r(z,t)\|_{\HH}^2\nu(dz)dt\Big)^{p}<\infty.$$
Especially, we denote $\cL_{\FF}^2((0,T)\times U_0; \HH)\equiv \cL_{\FF}^{2,1}((0,T)\times U_0;\HH)$.
   \item $\cL_\cF^{\nu, \infty} ([0,T]\times U_0; \HH)$: the space of $\cF_t$-predictable  bounded processes  $r=\{r(t,\omega,z),z\in U_0\}$  satisfying
   $$\int_{U_0}\|r(t,\omega, z)\|_{\HH}^p\nu(dz)\leq K,\  \forall \ p\geq 2.$$
\end{itemize}

We present the existence  and uniqueness for the solutions of   SDE and  BSDE  driven by the one-dimensional Brownian motion $W(t)$ and the Poisson random  measure ${N}(dt,dz)$ by the following two theorems, which are similar to Lemma 4.1 and Theorem 3.1 in Shen and Siu \cite{ST2013},  but will not be proved here. 
\begin{theorem}
  Consider the  MF-SDE with jumps 
  \begin{eqnarray}\label{GSDE}
 X(t)&=&a+\int_0^tb(s, X(s),\bar{X}(s), u(s), \bar{u}(s))ds+\int_0^t\sigma (s, X(s),\bar{X}(s), u(s), \bar{u}(s))dW_s\nonumber\\ 
   &&+\int_0^t\int_{U_0}\eta(s, X(s-),\bar{X}(s-), u(s-), \bar{u}(s-), z)\tilde{N}(ds,dz) 
\end{eqnarray}
where $a$ is an $\cF_0$-measurable random variable and $b, \si:[0,T]\times\Omega\times\RR^n\times\RR^n\times \RR^m\times\RR^m\to \RR^n$, $\eta:[0,T]\times\Omega\times\RR^n\times\RR^n\times \RR^m\times\RR^m\times U_0\to \RR^n$ are given  $\cF_t$-adapted processes satisfying the following assumptions 
\begin{enumerate}
    \item [(H1)]\label{1} There exists a constant $K>0$ such that for any $x, x'\in\RR^n, u,u'\in\RR^m,z\in U_0$ and a.s. $\omega\in \Omega$,
    \begin{eqnarray*}
    &&|b(t,x,x', u,u')|^2+|\si(t,x,x', u,u')|^2+\int_{U_0}|\eta(t,x,x', u,u',z)|^2\nu(dz)\\
    &\leq& K(1+|x|^2+|x'|^2+|u|^2+|u'|^2);\end{eqnarray*}
    \item [(H2)]  $b,\si$ and $\eta$  are uniformly Lipschitz continuous w.r.t. $x, x'$, i.e. there exist a constant $K>0$ such that for  all $u\in\cU, z\in U_0$,  $x,\ y,\ x',\  y'\in\RR^n$ and a,s, $\omega\in\Omega$, 
    \begin{eqnarray*}
&&|b(t,x, x', u,  u')-b(t,y, y', u, u')|^2 +|\si(t,x, x', u, u')-\si(t,y, y', u,  u')|^2 \nonumber\\
&&\ \ \ +\int_{U_0}|\eta(t,x, x', u,  u',z)-\eta(t,y, y', u, u',z)|^2\nu(dz)\leq K\big(|x-y|^2+|x'-y'|^2\big).
    \end{eqnarray*}
\end{enumerate}
Then, for any $u\in L_{\FF}^2(0,T; \RR^m)$,  the equation (\ref{GSDE}) admits a unique solution $X\in S^2_\FF(0,T; \RR^n)$. 
\end{theorem}

\begin{theorem}\label{2}
  Consider the  MF-BSDE with jumps 
  \begin{eqnarray}\label{GBSDE}
 Y(t)&=&\xi+\int_t^T f(s, Y(s),\bar{Y}(s), u(s), \bar{u}(s), Z(s),\bar{Z}(s), K(s), \bar{K}(s))ds\nonumber\\
&&-\int_t^TZ(s)dW(s)-\int_t^T\int_{U_0}K(s-,z)\tilde{N}(ds,dz), \  
\end{eqnarray}
where $\xi\in\RR^n$ is an $\cF_T$-measurable random variable and $f    
 :[0,T]\times\Omega\times\RR^n\times\RR^n\times \RR^m\times\RR^m\times\RR^{n}\times\RR^{n}\times\cL^2(U_0; \RR^n)\times\cL^2(U_0; \RR^n)\to\RR^n$  is a mapping  satisfying the following assumptions 
\begin{enumerate}
    \item [(H3)] $f(\cdot,0,0,\cdot,\cdot, 0,0,0,0)\in L_\FF^2(0,T; \RR^n)$ and $\xi\in L^2_{\cF_T}(\Omega; \RR^n)$. 
    \item [(H4)]  $f$  is  uniformly Lipschitz continuous w.r.t. $(y, z,k)$ and $(y', z', k')$, i.e. there exist a constant $K>0$ such that for  all $u\in \cU, z\in U_0$,  $\forall y,\ z,\ k,  \ y', z', k' \in\RR^n\times\RR^{n}\times\RR^n\times\RR^n\times\RR^{n}\times\RR^n$, and 
 $\forall y_1,\ z_1,\ k_1,  \ y_1', z_1', k_1' \in\RR^n\times\RR^{n\times n}\times\RR^n\times\RR^n\times\RR^{n}\times\RR^n$  , 
    \begin{eqnarray*}
&&|f(t, y,  y', u, u',  z,   z',  k, k' )-f(t, y_1, y_1', u, u', z_1,  z_1',  k_1,    k_1')|^2 \nonumber\\
&&\ \  \leq K\Big(|y-y_1|^2+|y'-y_1'|^2+|z-z_1|^2+|z'-z_1'|^2\\
&&\ \ \ +\int_{U_0}\(|k-k_1|^2+|k'-k_1'|^2\)\nu(dz)\Big), \ a,s, \ \omega\in\Omega 
    \end{eqnarray*}
\end{enumerate}
Then,  for any given $u\in L^2_{\FF}(0,T; \RR^m)$ , the equation (\ref{GBSDE})   has a unique solution $(Y, Z, K) \in S^2_\FF(0,T; \RR^n)\times L^2_\FF(0,T; \RR^{n})\times\cL_\FF^2([0,T]\times U_0; \RR^n)$.
\end{theorem}


Theorems 2.1 and 2.2  will play an important role in our paper for the existence and uniqueness of  solutions to some SDE and BSDE.

\section{  Main results } \label{sec3}
\setcounter{equation}{0}
\renewcommand{\theequation}{\thesection.\arabic{equation}}

In this section we present our main results, which give the optimal control  of Problem (MF-SLQ) by decomposing the problem into a constrained SLQ control problem without mean field and a deterministic  control problem with operator coefficients.

Throughout this paper, we make the following assumptions on the coefficients. 
\begin{itemize}
\item[(A1)]\label{A1} :  $A,\  A_1,\ C, \  C_1\in L_\FF^\infty (0,T; \RR^{n\times n})$, and  $B,\ \ B_1, \ D,\  D_1\in L_\FF^\infty (0,T; \RR^{n\times m})$, 
\item [(A2)]\label{A2} : $\al, \  \al_1\in \cL_\FF^{\nu,  \infty}([0,T]\times U_0; \RR^{n\times n}) $,  and $\be, \  \be_1\in \cL_\FF^{\nu, \infty}([0,T]\times U_0; \RR^{n\times m})$
\item[(A3)]\label{A3} : $Q, \  Q_1\in L_\FF^\infty(0,T; \cS_+^{n})$, $R, \  R_1\in L_\FF^\infty(0,T; \cS_+^{m})$ and $G\in L_{\cF_T}^\infty(\Omega, \cS_+^n)$. Moreover, the  weighting matrix processes $R$ are a.s. a.e. uniformly positive, i.e. there exist a constant $\delta>0$ such that 
\[R(t)\geq \delta I_m, \ \ a.e.\ t\in[0,T], \ \ a.s., \]
where $I_m$ is the $m\times m$ identity matrix.
\end{itemize}

Under the assumptions (A1)-(A2), it is obvious from Theorem 2.1 that for any given $u\in \cU_{ad}$ the state equation (\ref{sate}) has a unique solution $X\in S_{\FF}^2(0,T; \RR^n)$.  

 We  define  the Hamiltonian function $\cH:[0,T]\times\RR^n\times\RR^n\times\RR^m\times\RR^m\times\RR^n\times\RR^{n}\times\cL^2(U_0; \RR^n)\to \RR $ by
 \begin{eqnarray*}
    \cH(t, x,\bar{x}, u,\bar{u}, y,z, k )&:=&
    \big[A(t)x+A_1(t) \bar{x}+B(t)u+B_1(t)\bar{u}\big]^\top y\\
    &&+\tr\Big[\big[C(t)x+C_1(t)\bar{x}+D(t)u(t)+D_1(t)\bar{x}] \big]^\top z\Big]\\
&&+\int_{U_0}\big[\al(t,z)x+\al_1(t,z)\bar{x}+\be(t,z)u+\be_1(t)\bar{u}\big]^\top k(z)\nu(dz)\\
&& +\frac{1}{2}\big(x^\top Q(t)x+  \bar{x}^\top Q_1(t)\bar{x}
 +u^\top R(t)u+\bar{u}^\top R_1(t)\bar{u}\big).
 \end{eqnarray*}
Let $(u, X)$ be an admissible pair.  The corresponding adjoint BSDE of the stochastic system  (\ref{sate}) is defined by 
\begin{equation}\label{adjion’}
\left\{\begin{array}{ccl}
dY(t)&=&-\cH_x(t,X(t),\bar{X}(t), u(t), \bar{u}(t), Y(t), Z(t), K(t))dt\\
&&\ \ - \EE\big[\cH_{\bar{x}}(t,X(t),\bar{X}(t), u(t), \bar{u}(t), Y(t), Z(t), K(t))\big]dt\\
&&\ \ +Z(t)dW(t)+\int_{U_0}K(t,z)\tilde{N}(dt,dz)\\
Y(T)&=&GX(T), \ t\in[0,T].
\end{array}\right.
\end{equation}
Note that under assumptions (A1)-(A2), Theorem \ref{2} implies that for any given $u\in \cU_{ad}$, the equation (\ref{adjion’}) admits a unique solution 
\[(Y,\ Z,\  K) \in S^2_\FF(0,T; \RR^n)\times L^2_\FF(0,T; \RR^{n})\times\cL_\FF^2([0,T]\times U_0; \RR^n).\]

Now we give two  preliminary results of  the optimal control
problem (MF-SLQ).
\begin{theorem}\label{UST}
    Let (A1)-(A3) hold. Then Problem (MF-SLQ) admits a unique optimal control $u^*\in\cU_{ad}.$ 
\end{theorem}

\begin{theorem}\label{pMSP}
 Let (A1)-(A3) hold. 
 A control $u^*$  is optimal for  Problem (M-SLQ) if and only if the solution $(X^*, Y^*,Z^*, K^*)$ to the FBSDE
\begin{equation}\label{FBSDE11}
\left\{\begin{array}{ccl}
dX^*(t)&=&\big(AX^*+A_1\bar{X}^*+Bu^*+B_1\bar{u}^*\big)dt +\big(CX^*+C_1\bar{X}^*+Du^*+D_1\bar{u}^* \big)dW(t)\\
&&\ \ +\int_{U_0}\big(\al(z)X^*+\al_1(z)\bar{X}^*+\be(z)u^*+\be_1\bar{u}^*\big)\tilde{N}(dt,dz)\\

dY^*(t)&=&-\Big(A^\top Y^*+C^\top Z^*+\int_{U_0}\al^\top(z)K^*(z)\nu(dz)+QX^*\Big)dt\\
&&\ \ -\Big(\EE\Big[A_1^\top Y^*+C_1^\top Z^*+\int_{U_0}\al_1^\top(z)K^*(z)\nu(dz)\Big]+\EE[Q_1]\bar{X}^*\Big)dt\\
&&\ \ +Z^*dW(t)+\int_{U_0}K^*(z)\tilde{N}(dt,dz)\\
X^*(0)&=&x,\ Y(T)=GX^*(T), \ t\in[0,T],
\end{array}\right.
\end{equation}
satisfies the stationary condition
\begin{eqnarray}\label{SCUB}
\cH_u(t, X^*, \bar{X}^*, u, \bar{u}^*, Y^*, Z^*, K^*)+\EE[\cH_{\bar{u}}(t, X^*, \bar{X}^*, u^*, \bar{u}^*, Y^*, Z^*, K^*)]=0,
\end{eqnarray}
for $a.e. \ t\in[0,T], \ a.s.$
\end{theorem}

In (\ref{pMSP})   above, we have omitted the time parameters in the drift term because the expression is too complicated.  We will continue to do so in the rest of this article when it is not confusing.  Furthermore,  for any  given $u^*\in\cU_{ad}$,   Theorems  \ref{GSDE} and   \ref{2}  in Section 2 imply that  the FBSDE (\ref{FBSDE11}) admits a unique solution 
\[(X^*, Y^*,\ Z^*,\  K^*) \in (S^2_\FF(0,T; \RR^n))^2\times L^2_\FF(0,T; \RR^{n})\times\cL_\FF^2([0,T]\times U_0; \RR^n).\]

Because of the terms such as $\EE\big[A_1^\top Y^*]\neq\EE[A_1^\top ] \EE[Y^*]$ in the adjoint equation, it is challenging to decouple the optimal system (\ref{FBSDE11}-\ref{SCUB}). We will adapt the  ELM method 
 introduced in \cite{XW2024} to the current setup to overcome this difficulty.  

  We consider  the constraints that the state process $X(t)$ and the admissible control $u(t)$ satisfy 
 $\EE[X(t)]=a(t)$ and $\EE[u(t)]=b(t)$, $\forall t\in[0,T]$, where $a,\;b$ are two fixed deterministic functions. 
 We begin by rewriting Problem  (MF-SLQ) as 
\begin{eqnarray*}
   \inf_{u\in\cU_{ad}}J(u)=\inf_{(a, b)\in L^2(0,T; \RR^{n+m})}\inf_{ u\in\cU_{ad}}\big\{J(u): \EE[X]=a, \EE[u]=b\big\}. 
\end{eqnarray*}

 First,  we consider a control problem with the constraints, which  called the Problem 1. In this case, we denote the cost functional as $J_{a, b}$. Namely, 
 \begin{eqnarray}
  J_{a,b}(u) &=&\EE\bigg\{\big<GX(T), X(T)\big>
+\int_0^T \Big(\big<Q(t)X(t), X(t)\big> +  \big<Q_1(t)a(t), a(t)\big>\nonumber\\
&&\ \ \ +\big<R(t)u(t), u(t)\big>+\big<R_1(t)b(t), b(t)\big> \Big)dt\bigg\},
 \end{eqnarray}
 where $X$ satisfies equation (\ref{sate}) with the constraints $\EE[X]=a$ and $\EE[u]=b$.   Then this constrained problem is formulated as follows:
 
 $\bf{Problem}$ $\bm{1.}$ 
 For a given initial state $x\in\RR^n$,  find a control $u^*_{a, b}\in \cU_{ad}$ such that 
 \begin{equation}\label{EJ2HH} 
 J_{a, b}(u^*_{a,b})=\inf_{u\in\cU_{ad}}\left\{J_{a, b}(u)\Big|\EE[X^u]=a, \EE[u]=b\right\}.
 \end{equation}
 
\begin{lemma}\label{P1U1}
For  $a\in L^2(0,T; \RR^n)$  and  $b\in L^2(0,T; \RR^m)$  fixed, there exists a unique $u^*_{a, b}\in \cU_{ad}$  such that $\EE[X^{u^*_{a, b}}]=a$,  $\EE[u^*_{a, b}]=b$ and 
\[J_{a, b}(u^*_{a, b})=\inf_{u\in\cU_{ad}}\big\{J_{a, b}(u):\EE[X^u]=a, \EE[u]=b\big\}.\]
 \end{lemma}

To solve  Problem 1,  we relax the constraint by introducing the ELMs $\la$ and $\ga$, which are  deterministic functions. Specifically, we consider a control problem with state equation:  
\begin{equation}\label{EXcsate}
\left\{\begin{array}{ccl}
dX(t)&=&\big(AX+A_1a+Bu+B_1b\big)dt+\big(CX+C_1a+Du+D_1b \big)dW(t)\\
&&\ \ +\int_{U_0}\big(\al(z)X+\al_1(z)a+\be(z)u+\be_1(z)b\big)\tilde{N}(dt,dz)\\
X(0)&=&x, \ t\in[0,T],
\end{array}\right.
\end{equation}
 and  the objective functional  given by 
 \begin {eqnarray}\label{Ccost1}
\hat{J}_{a,b}(u , \la , \ga ) &=&J_{a, b}(u)+2\big<\la, \EE[X]-a\big>_{L^2}+2\big<\ga, \EE[u]-b\big>_{L^2}.
\end {eqnarray}
  where $\lambda(t)\in L  ^2(0,T; \RR^n)$ and $\gamma(t)\in L^2(0,T; \RR^n)$ are deterministic functions.

 Note that  $\cU_{ad}$ and  $L^2(0,T;\RR^{n+m})$ are convex sets, the function   $\hat{J}_{a, b}(u,\la, \ga)$ is directionally differentiable  
 with respect to $(u,\la, \ga)$ and  for any $(u, \la, \ga)\in \cU_{ad}\times L^2(0,T; \RR^{n})\times L^2(0,T; \RR^{m})$ the functions $\hat{J}_{a, b}(\cdot,\la, \ga)$ and $\hat{J}_{a, b}(u, \cdot, \cdot)$  are convex and concave, respectively.   According to  Propositions 2.157  and  2.156 in \cite{BS2013} there is no duality gap between Problem 1 and its dual.  Then, Problem 1 will be solved through two subproblems. Our first  is  fixing $\la\in L^2(0,T; \RR^n)$ and $\ga  \in  L^2(0,T; \RR^m) $   to find $u_{a, b, \la, \ga}\in \cU_{ad}$ such that the objective functional (\ref{Ccost1}) is minimized.  
 Once we have found the  optimal control $u_{a, b, \la, \ga}$, the second   subproblem is to find  the optimal ELMs $\la_{a, b}$  and $\ga_{a, b}$ such that the cost functional $\hat{J}_{a, b}( u_{a, b, \la, \ga},\la,\ga)$ is maximised.   For a more precise expression, we write $\hat{J}_{a,b}( u, \la, \ga)$ as 
 $J_{a,b, \la, \ga}(u)$.  
 
$\bf{Problem}$ $\bm{(Sub1).}$ 
 For a given initial state $x\in\RR^n$, $a, \la\in L^2(0,T; \RR^n)$ and $b, \ga\in L^2(0,T; \RR^m)$, find a control $u_{a, b, \la, \ga}\in \cU_{ad}$ such that 
 \begin{equation}\label{J2HH} 
 J_{a,b,\la, \ga}(u_{a,b, \la, \ga})=\inf_{u\in\cU_{ad}}J_{a,b,\la, \ga}( u).
 \end{equation}

 \begin{lemma}\label{UN2} 
Assume that (A1)-(A3) hold. Then,  for any $x\in\RR^n$, $a, \la\in L^2(0,T; \RR^n)$ and $b, \ga\in L^2(0,T; \RR^m)$ fixed, 
  Problem (Sub1)  admits a unique optimal control $u_{a, b, \la, \ga}\in \cU_{ad}$.
 \end{lemma} 

  \begin{lemma}\label{OPSY2} 
Assume that (A1)-(A3) hold. Then,  for any $x\in\RR^n$, $a, \la\in L^2(0,T; \RR^n)$ and $b, \ga\in L^2(0,T; \RR^m)$ fixed, 
the  control $u_{a, b, \la, \ga}$  is optimal if and only if  the  solution $(X_{a, b, \la, \ga}, Y_{a, b, \la, \ga}, Z_{a, b, \la, \ga})$
   to the  FBSDE 
\begin{equation}\label{Cadjion}
\left\{\begin{array}{ccl}
dX_{a,b, \la, \ga}&=&\big(AX_{a,b, \la, \ga}+A_1a+Bu_{a,b, \la, \ga}+B_1b\big)dt \\
&&+\big(CX_{a,b, \la, \ga}+C_1a+Du_{a,b, \la, \ga}+D_1b \big)dW(t)\\
&&\ \ +\int_{U_0}\big(\al(z)X_{a,b, \la, \ga}+\al_1(z)a+\be(z)u_{a,b, \la, \ga}+\be_1(z)b\big)\tilde{N}(dt,dz)\\
dY_{a,b, \la, \ga}(t)&=&-\big(A^\top Y_{a,b, \la, \ga}+C^\top Z_{a,b, \la, \ga}+\int_{U_0}\al^\top(z)K_{a,b, \la, \ga}(z)\nu(dz)+Q X_{a,b, \la, \ga}+\la \big)dt\\
&&\ \ +Z_{a,b, \la, \ga}(t)dW(t)+\int_{U_0}K_{a,b, \la, \ga}(z)\tilde{N}(dt,dz)\\
X_{a,b, \la, \ga}(0) &=&x, \ Y_{a,b, \la, \ga}(T)=GX_{a,b, \la, \ga}(T), \ t\in[0,T].
\end{array}\right.
\end{equation} 
 satisfies the  the stationary condition 
\begin{equation}\label{HST}
B^\top Y_{a,b, \la, \ga}+D^\top Z_{a,b, \la, \ga}+\int_{U_0}\be^\top K_{a,b, \la, \ga}(z)\nu(dz)+\ga+Ru_{a,b, \la, \ga}=0.
\end{equation}
 \end{lemma} 

The system (\ref{Cadjion}), together with the stationary condition (\ref{HST}), is referred to as the optimality system for Problem (Sub1).  Substituting (\ref{HST}) into  (\ref{Cadjion}), the optimality system becomes a coupled FBSDE as follows:
\begin{equation}\label{OMCadjion} 
\left\{\begin{array}{ccl}
dX_{a,b, \la, \ga}&=&\Big(AX_{a,b, \la, \ga}+A_1a+B_1b\\
&&\ \ \  -BR^{-1}\(B^\top Y_{a,b, \la, \ga}+D^\top Z_{a,b, \la, \ga}
+\int_{U_0}\be^\top K_{a,b, \la, \ga}(z)\nu(dz)+\ga\)\Big)dt \\
&&+\Big(CX_{a,b, \la, \ga}+C_1a+D_1b\\
&&\ \ \ -DR^{-1}\big(B^\top Y_{a,b, \la, \ga}+D^\top Z_{a,b, \la, \ga}+\int_{U_0}\be^\top K_{a,b, \la, \ga}(z)\nu(dz)+\ga\big)\Big)dW(t)\\
&& +\int_{U_0}\Big(\al(z)X_{a,b, \la, \ga}+\al_1(z)a +\be_1(z)b\\
&&\ \ \ \ -\be(z)R^{-1}\big(B^\top Y_{a,b, \la, \ga}+D^\top Z_{a,b, \la, \ga}+\int_{U_0}\be^\top K_{a,b, \la, \ga}(\tilde z)\nu(d\tilde z)+\ga \big)\Big)\tilde{N}(dt,dz)\\
dY_{a,b, \la, \ga}(t)&=&-\Big(A^\top Y_{a,b, \la, \ga}+C^\top Z_{a,b, \la, \ga}+\int_{U_0}\al^\top(z)K_{a,b, \la, \ga}(z)\nu(dz)+Q X_{a,b, \la, \ga}  +\la \Big)dt\\
&&\ \ +Z_{a,b, \la, \ga}(t)dW(t)+\int_{U_0}K_{a,b, \la, \ga}(z)\tilde{N}(dt,dz)\\
X_{a,b, \la, \ga}(0) &=&x, \ Y_{a,b, \la, \ga}(T)=GX_{a,b, \la, \ga}(T), \ t\in[0,T].
\end{array}\right.
\end{equation}

Combining  Lemma \ref{UN2} with Lemma \ref{OPSY2},  we get the unique solvability of FBSDE (\ref{OMCadjion}).
\begin{lemma}\label{LOP1}
  Let (A1)-(A3) hold.  Then,  for any $x\in\RR^n$, $a, \la\in L^2(0,T; \RR^n)$ and $b, \ga\in L^2(0,T; \RR^m)$  fixed,   the coupled system  (\ref{OMCadjion})  admits a unique adapted solution
$$(X_{a, b, \la, \ga}, Y_{a, b, \la, \ga}, Z_{a, b, \la, \ga}, K_{a, b, \la, \ga})\in (S_{\FF}^2(0,T; \RR^n))^2\times L_{\FF}^2(0,T; \RR^n)\times \cL_{\FF}^2([0,T]\times U_0; \RR^n).$$ 
\end{lemma}

In order to  obtain the explicit expression  of  the optimal control $u_{a, b, \la, \ga}$ for Problem (Sub1), we need to  decouple the optimality system (\ref{OMCadjion}). To simplify the notation, we denote
 \begin{eqnarray*}  
\Theta&=&-\Big(D^\top P  D +\int_{U_0}\Big(\be^\top(z) \big(P+\Gamma(z)) \be(z)\Big)\nu(dz)+R\Big);\ \nonumber\\
N&=&B^\top  P +D^\top \Lambda +D^\top  P  C  +\int_{U_0}\Big(\be^\top(z)\Gamma(z)+\be^\top(z)\big(P +\Gamma(z)\big) \al(z)\Big)\nu(dz);\nonumber\\
   H&=& C^\top P D+P B+\Lambda D+\int_{U_0}\Big(\al^\top(z) \big(P+\Gamma(z)\big)\be(z)+\Gamma(z)\be(z)\Big)\nu(dz);\nonumber\\
\hat{A}&=&A+B\Theta^{-1}N;\ \  \hat{A}_1=A_1+B\Theta^{-1}D^\top P C_1+\int_{U_0}B\Theta^{-1}\be^\top(z) \big(P+\Gamma\big) \al_1(z)\nu(dz);\nonumber\\
\hat{B}_1&=&B_1+B\Theta^{-1}D^\top P D_1+\int_{U_0}B\Theta^{-1}\be^\top(z)\big(P+\Gamma\big)  \be_1(z)\nu(dz);\nonumber\\
\hat{C}&=&C+D\Theta^{-1}N; \ \hat{C}_1=C_1+D\Theta^{-1}D^\top P C_1+\int_{U_0} D\Theta^{-1}\be^\top(z) \big(P+\Gamma\big)  \al_1(z)\nu(dz);\nonumber\\
\hat{D}_1&=& D_1+D\Theta^{-1}D^\top P D_1+\int_{U_0}D\Theta^{-1}\be^\top(z) \big(P+\Gamma\big)  \be_1(z)\nu(dz);\nonumber\\
\hat{\al}(z)&=&\al(z)+\be(z)\Theta^{-1}N;\nonumber\\
\hat{\al}_1(z)&=&\al_1(z)+\be(z)\Theta^{-1}D^\top P C_1
+\int_{U_0}\be(\tilde z)\Theta^{-1}\be^\top(\tilde z) \big(P+\Gamma\big)  \al_1(\tilde z)\nu(d\tilde z);\nonumber\\
\hat{\be}_1(z)&=&\be_1(z) +\be(z)\Theta^{-1}D^\top \big(P+\Gamma\big)  D_1
+\int_{U_0}\be(\tilde z)\Theta^{-1}\be^\top(\tilde z) P \be_1(\tilde z)\nu(d\tilde z);\nonumber\\
\hat{M}&=&A^\top + H\Theta^{-1}B^{\top}; \ \  \hat{N}=C^\top+H\Theta^{-1}D^\top\,; \  \hat{K}(z) =\al(z)+H\Theta^{-1}\be^{\top}(z);\nonumber\\
\hat{L}&=& C^\top P C_1+P A_1+\Lambda C_1+\int_{U_0}\big(\al^\top(z) P\al_1(z)+\al^\top(z) \Gamma(z)\al_1(z)+\Gamma(z)\al_1(z)\big)\nu(dz)\nonumber\\
&&\ \ \ 
+H\Theta^{-1}D^\top P C_1+\int_{U_0}H\Theta^{-1}\big(\be^\top P \al_1+\be^\top \Gamma(z)\al_1\big)\nu(dz);\nonumber\\
\Hat{Q}&=&C^\top P D_1+P B_1+\Lambda D_1+\int_{U_0}\big(\al^\top(z) P\be_1(z)+\al^\top(z) \Gamma(z)\be_1(z)+\Gamma(z)\be_1(z)\big)\nu(dz)\nonumber\\
&&\ \ \ 
+H\Theta^{-1}D^\top P D_1+\int_{U_0}H\Theta^{-1}\big(\be^\top(z) P \be_1(z)+\be^\top(z) \Gamma(z) \be_1(z)\big)\nu(dz).
 \end{eqnarray*}

  First,  we introduce the following decoupling equations for the FBSDE (\ref{OMCadjion}) which will be discussed in detail in Section \ref{sec51}. We consider the SREJ 
 \begin{equation}\label{RP}
\left\{\begin{array}{ccl}
dP(t)&=&-\Big\{A^\top P+P A+C^\top\Lambda+\Lambda C+C^\top P  C+Q+N^\top \Theta^{-1} N\\
&&\ \ \ +\int_{U_0}\Big(\al^\top(z)\Gamma(z)+\Gamma(z)\al(z)+\al^\top(z) \big(P+\Gamma(z)\big) \al(z)\Big)\nu(dz)\Big\}dt
\\
&&+\Lambda (t)dW(t)+\int_{U_0} \Gamma(t-,z)\tilde{N}(dt,dz), \qquad t\in[0,T]\\
P(T)&=&G,
\end{array}\right.
\end{equation}
 and  the  linear BSDE with jumps
\begin{equation}\label{Lphi} 
\left\{\begin{array}{ccl}
d\phi(t)&=&-\bigg\{\hat{M}\phi+\hat{N}\psi +\int_{U_0}\hat{K}(z)\theta(z)\nu(dz)+\hat{L}a+\hat{Q}b+\lambda+H\Theta^{-1}\ga
\bigg\}dt\\
&&+\psi(t) dW(t)+\int_{U_0}\theta(t-,z)\tilde{N}(dt,dz)\\
\phi(T)&=&0, \ t\in[0,T].
\end{array}\right.
\end{equation}

Note that the  SREJ (\ref{RP})  is a highly  nonlinear BSDE with jumps, the general theorem of BSDE (see theorem \ref{2})  cannot be applied to  guarantee the existence and uniqueness of its solution. However,  this equation  has  been  introduced  by Meng \cite{Meng2014}  in order to decouple the FBSDE (\ref{OMCadjion}) for  the homogeneous  case (i.e. $a=b=\la=\ga=0$), and the existence and uniqueness of solution to SREJ (\ref{RP}) have also  been obtained by Zhang et al. \cite{ZDM2020}. The following Lemma  can be derived directly from Theorem 4.1 and Theorem 5.1 in the work of Zhang et al \cite{ZDM2020}.

 \begin{lemma}\label{BSRE}
  Assume  that (A1)-(A3) hold. Then the stochastic Riccati equation (\ref{RP})  admits a unique solution $(P, \Lambda, \Gamma)\in S_\FF^\infty(0,T; \cS^n_+)\times L_\FF^2(0,T; \cS^n)\times \cL_\FF^2\big([0,T]\times U_0;  \cS^n\big)$ 
  such that for some $c>0$, 
\[D^\top P D+\int_{U_0}\be^\top(z) \big(P+\Gamma(z)\big)\be(z)\nu(dz)+R \geq cI_{n}, \ a.e \ on\ t\in[0,T], \ a.s.\]

Moreover, there is a  constant $K$ such that the following estimate holds:
\[\EE\Big(\int_0^T|\Lambda(t)|^2dt+\int_0^T\int_{U_0}|\Gamma(t,z)|^2\nu(dz)dt\Big)\leq K.\]
\end{lemma}
\bigskip

Although (\ref{Lphi})  is a linear BSDE, its coefficients are unbounded since  $\Lambda$   and $\Gamma(z)$ are so. If there is no jump term, i.e. $\Gamma(z)=0$,   for the one-dimensional case, the existence and  uniqueness of the solution to (\ref{Lphi}) has been obtained by Hu et al \cite{HSX} using the BMO-martingale and Girsanov Theorem.   Since the  Girsanov Theorem of jump-processes needs to require $\hat{K}(z)>-1$, the method in \cite{HSX} is not applicable. In this  article, we decouple (\ref{OMCadjion}) for   non-homogeneous case with  L\'evy driving noises of any dimensions without any additional conditions besides (A1)-(A3) and prove the existence and uniqueness of  solution to the   linear BSDE  (\ref{Lphi}) in Section \ref{sec51}. 

\begin{lemma}\label{LBSDER}
Let (A1)-(A3) hold. Then  the  decoupling equation (\ref{Lphi}) for the FBSDE (\ref{OMCadjion}) admits a unique solution $(\phi, \psi, \theta)\in S_{\FF}^2(0,T; \RR^n)\times L_{\FF}^{2, p}(0,T; \RR^n)\times \cL_{\FF}^{2, p}([0,T]\times U_0; \RR^n), \forall p\in(1/2,1). $
\end{lemma}

Combining the above results, we can obtain the following theorem, which gives the  state feedback representation of the
 optimal control $u_{a, b, \la, \ga}$  of Problem (Sub1). 
 
 \begin{theorem}\label{TH2'}
Let  (A1)-(A3) hold.   Then,  the  optimal control $u_{a,b, \la, \ga}(t)$  of  Problem (Sub1) has a feedback  representation of the state $X_{a,b, \la, \ga}(t)$, 
\begin{equation}\label{OPu'''}
u_{a,b, \la, \ga}=\Theta^{-1}N X_{a,b, \la, \ga}
+\Theta^{-1}M, 
\end{equation}
 where 
 \begin{eqnarray}\label{MPAO1}
M&=&B^\top \phi +D^\top \psi+D^\top P C_1a+D^\top P D_1 b+\ga\nonumber\\
&&\ \ \ \ +\int_{U_0}\Big(\be^\top \theta+\be^\top \big(P+\Gamma(z)\big) \al_1 a+\be^\top\big(P+\Gamma(z)\big) \be_1 b \Big)\nu(dz), 
\end{eqnarray}
$$(P, \Lambda, \Gamma)\in S_\FF^2(0,T, S^n)\times L_\FF^2(0,T, S^n)\times \cL_\FF^2(0,T, S^n),$$
and  $$(\phi, \psi, \theta)\in S_\FF^2(0,T, \RR^n)\times L_\FF^{2, p}(0,T, \RR^n)\times \cL_\FF^{2, p}(0,T, \RR^n),\quad \big(p\in(1/2,1)\big)$$ are the unique solutions of  SREJ (\ref{RP})  and linear  BSDE  (\ref{Lphi}), respectively. 
 
 Furthermore, the optimal state $X_{a,b, \la, \ga}$ is the unique solution of the following  SDE
 \begin{equation}\label{X0}
\left\{\begin{array}{ccl}
dX_{a,b, \la, \ga}(t)&=&\bigg(B\Theta^{-1}\Big(B^\top\phi+D^\top\psi
 +\int_{U_0}\be^\top(z) \theta(z)\nu(dz)\Big) \\
&&\ \ \ +\hat{A}X_{a,b, \la, \ga}+\hat{A}_1 a+\hat{B}_1 b+B\Theta^{-1}\gamma\bigg)dt\\
&&+\bigg(D\Theta^{-1}\Big(B^\top \phi +D^\top\psi
  +\int_{U_0}\be^\top(z)\theta(z)\nu(dz)\Big)
  \\
&&\ \ \ +\hat{C}X_{a,b, \la, \ga}+\hat{C}_1 a+\hat{D}_1b+D\Theta^{-1}\gamma\bigg)dW(t)\\
&& +\int_{U_0}\bigg(\be(z)\Theta^{-1}\Big(B^\top \phi +D^\top \psi+\int_{U_0}\be^\top(\tilde z)\theta(\tilde z)\nu(d\tilde z)\Big)
\\
&&\ \ \ +\hat{\al}(z)X_{a,b, \la, \ga}+\hat{\al}_1 a+\hat{\be}_1 b+\be\Theta^{-1}\gamma\bigg)\tilde{N}(dt,dz),\\
X(0)&=&x, \ \ \ t\in[0,T].
\end{array}\right.
\end{equation}
\end{theorem}

\bigskip
 The existence and  uniqueness of the equations (\ref{Lphi}-\ref{X0})  imply that for any $x\in  \RR^n$ $a, \lambda\in L^2(0,T; \RR^n)$ and $b, \gamma\in L^2(0,T; \RR^m)$, 
 the equations (\ref{Lphi}) and (\ref{X0}) have a unique solution 
 $$\big(X_{( a, b, \lambda, \gamma)}, \phi_{( a, b, \lambda, \ga)}, \psi_{ (a, b, \lambda, \ga)}, \theta_{ (a, b, \lambda, \ga)}\big)$$ 
 taking values in  $S_\FF^2(0,T; \RR^n) \times S_\FF^2(0,T; \RR^n)\times L_\FF^{2, p}(0,T; \RR^n)\times\cL_\FF^{2, p}([0,T]\times U_0; \RR^n)$ with $p\in(1/2, 1)$.
  We define linear operators $\cL_i$, $i=0,1,\cdots,4$, from $\RR^n$, $L^2(0,T; \RR^n)$, $L^2(0,T; \RR^m)$,  $L^2(0,T; \RR^n)$, $L^2(0,T; \RR^m)$, respectively, to $L^2(0,T; \RR^n)$  such that 
\begin{eqnarray}\label{BSX1}
   \EE[X_{a, b, \la, \ga}(t)] =(\cL_0 x)(t)+(\cL_1 a)(t)+(\cL_2 b)(t)+ (\cL_3 \la)(t)+(\cL_4 \ga)(t).
\end{eqnarray}

 We  also introduce the following processes
\begin{eqnarray}\label{Ma}
M_{a}&=& B^\top \phi_{(a, 0,0, 0)}+D^\top \psi_{(a, 0,0, 0)}+D^\top P C_1a\nonumber\\
&&+\int_{U_0}\be^\top(z) \(\theta_{(a, 0,0, 0)}(z)+ P \al_1(z) a\)\nu(dz),
\end{eqnarray}
\begin{eqnarray}\label{Mb}
 M_b&=& B^\top \phi_{(0, b,0, 0)}+D^\top \psi_{(0, b,0, 0)}\nonumber\\
 &&+D^\top PD_1 b+\int_{U_0}\be^\top(z) \(\theta_{(0, b,0, 0)}(z)+P \be_1(z) b\)\nu(dz),
\end{eqnarray}
\begin{eqnarray}\label{Mb}
 M_\lambda&=& B^\top \phi_{(0, 0,\lambda, 0)}+D^\top \psi_{(0, 0,\lambda, 0)}+\int_{U_0}\be^\top(z) \theta_{(0, 0,\lambda, 0)}(z)\nu(dz),
\end{eqnarray}
 and
 \begin{eqnarray}\label{ML}
 M_\gamma&=& B^\top \phi_{(0, 0,0, \ga)}+D^\top \psi_{(0, 0,0, \ga)}+\int_{U_0}\be^\top(z) \theta_{(0, 0,0, \ga)}(z)\nu(dz)+\gamma.
\end{eqnarray}
 
  We observe that in the above systems $(\ref{Ma})-(\ref{ML})$, $M_a, M_b, M_\lambda,  M_\gamma$ can be represented explicitly as linear functionals of $a, b,\lambda$ and $\gamma$.  Recall that $u_{a, b, \la, \ga}$ and $M$ are given in (\ref{OPu'''}) and (\ref{MPAO1}),  respectively.  In fact, there also exist linear operators $\tilde{\cL}_0, \tilde{\cL}_1, \tilde{\cL}_2,  \tilde{\cL}_3 $  and $\tilde{\cL}_4$  from $\RR^n$, $L^2(0,T; \RR^n)$, $L^2(0,T; \RR^m)$,  $L^2(0,T; \RR^n)$, $L^2(0,T; \RR^m)$, respectively, to $L^2(0,T; \RR^m)$  such that  
\begin{eqnarray}\label{BSU1}
  \EE[u_{a, b, \la, \ga}(t)]=(\tilde{\cL}_0x)(t)+(\tilde{\cL}_1a)(t)+(\tilde{\cL}_2b)(t)+(\tilde{\cL}_3\la)(t)+(\tilde{\cL}_4\ga)(t). 
\end{eqnarray}

We  define the function $\tilde J_{a,b}(\la, \ga):L^2([0,T],\RR^n)\times L^2([0,T],\RR^m)\to \RR$ as the minimum value of the cost function   $J_{a,b, \la, \ga}(u)$ over $u\in\cU_{ad}$: for $a\in L^2([0,T],\RR^n), \  b \in L^2([0,T],\RR^m)$, 
\begin{eqnarray}\label{DCost1'}
&&\tilde J_{a,b}(\la, \ga)=\inf_{u\in\cU_{ad}}J_{a, b, \lambda, \gamma}( u )\nonumber\\
 &=&\EE[\big<G X_{a,b, \la, \ga}(T), X_{a,b, \la, \ga}(T)\big>] +\EE\int_0^T\bigg\{\big<QX_{a,b, \la, \ga}, \ X_{a,b, \la, \ga}\big>+\big<Q_1a, \ a\big> +\big<R_1b, \ b\big>
\nonumber\\
 &&\   +2\big<\lambda, X_{a,b, \la, \ga}-a\big>+2\big<\ga, \ -b\big>-\big<R^{-1}\gamma, \ \gamma \big>+\big<BR^{-1}B^{\top}Y_{a,b, \la, \ga}, \ Y_{a,b, \la, \ga} \big>\nonumber\\
 &&\ \  +\big<D R^{-1} D^\top Z_{a,b, \la, \ga}, \ Z_{a,b, \la, \ga} \big> +\big<R^{-1}\int_{U_0}\be^\top(z)K_{a,b, \la, \ga}(z)\nu(dz), \ \int_{U_0}\be^\top(z)K_{a,b, \la, \ga}(z)\nu(dz) \big>\nonumber\\
 &&\ \ +2\big<D R^{-1}B^\top Y_{a,b, \la, \ga}, \  Z_{a,b, \la, \ga}\big>+2\big<R^{-1}B^\top Y_{a,b, \la, \ga}, \ \int_{U_0}\be^\top(z)K_{a,b, \la, \ga}(z)\nu(dz)\big>\nonumber\\
 &&\ \ +2\big<R^{-1}D^\top Z, \ \int_{U_0}\be^\top(z)K_{a,b, \la, \ga}(z)\nu(dz)\big> \bigg\}dt,
\end {eqnarray}
 where $(X_{a, b, \la, \ga}, Y_{a, b, \la, \ga}, Z_{a, b, \la, \ga}, K_{a, b, \la, \ga})$ is the  unique solution of system (\ref{OMCadjion}). We consider the following
 
$\bf{Problem}$ $\bm{(Sub2).}$ 
 For given initial state $x\in\RR^n$, $a\in L^2([0,T],\RR^n)$ and $b\in L^2([0,T],\RR^m)$, find  optimal 
ELMs $\la_{a, b}\in L^2(0,T; \RR^{n})$ and $ \ga_{a, b}\in L^2(0,T; \RR^{m})$ such that 
 \begin{equation}\label{V1HH} 
\tilde J_{a, b}(\la_{a, b}, \ga_{a, b})=\sup_{(\la, \ga)\in L^2(0,T; \RR^{n+m})}\tilde J_{a,b}(\la,\ga).
 \end{equation}
 Now we  also give the   maximum principle  of  the  Problem (Sub2).

\begin{lemma}\label{DMSP}
Assume that (A1)-(A3) hold.    Let $(\la_{a, b}, \ga_{a, b})$  be an optimal solution for Problem (Sub2).
Then  the  solution $(X_{a, b, \la_{a, b}, \ga_{a, b}}, Y_{a, b, \la_{a, b}, \ga_{a, b}}, \  Z_{a, b, \la_{a, b}, \ga_{a, b}}, \ K_{a, b, \la_{a, b}, \ga_{a, b}})$ to  the linear jump-diffusion FBSDE (\ref{OMCadjion}) 
 satisfies (with  $\la$ and $\ga$ replaced by $\la_{a, b}$ and $\ga_{a, b}$, respectively)  the stationary conditions
 \begin{equation}\label{DOS1}
\EE[u_{a, b, \la_{a, b}, \ga_{a, b}}]-b=0,
\end{equation}
and
 \begin{equation}\label{DOS2}
\EE[ X_{a, b, \la_{a, b}, \ga_{a, b}}]-a=0,
\end{equation}
where $(X_{a, b, \la_{a, b}, \ga_{a, b}},Y_{a, b, \la_{a, b}, \ga_{a, b}},Z_{a, b, \la_{a, b}, \ga_{a, b}},K_{a, b, \la_{a, b}, \ga_{a, b}})$ is the  unique solution of the  system (\ref{OMCadjion}) with $\la$ and  $\ga$  replaced by $\la_{a, b}$ and $\ga_{a, b}$, respectively.
\end{lemma}

\begin{theorem}\label{BUAB}
     For any  EMLs $\la_{a, b}$ and $\ga_{a, b}$ such that   $\EE[X_{a, b, \la_{a, b}, \ga_{a, b}}]=a$ and $\EE[u_{a, b, \la_{a, b}, \ga_{a, b}}]=b$, we have $u_{a, b, \la_{a, b}, \ga_{a, b}}=u^*_{a, b}$.  Furthermore,   the optimal control $u^*_{a,b}$ of Problem 1  is given by (\ref{OPu'''}) with the  optimal state $X^*_{a, b}$  satisfying (\ref{X0}) (with  $\la$ and  $\ga$  replaced by any solution of equations (\ref{DOS1}) and (\ref{DOS2})).
\end{theorem}

 Recall that some operators representations of $X_{a, b, \la, \ga}$ and $u_{a, b, \la, \ga}$ are given by (\ref{BSX1}) and (\ref{BSU1}). Thus,
\begin{eqnarray}\label{EXB}
  (\cL_0x)(t)+(\cL_1 a)(t)+(\cL_2 b)(t)+ (\cL_3 \la_{a, b})(t)+(\cL_4 \ga_{a, b})(t)=a(t).
\end{eqnarray}
 and  
 \begin{eqnarray}\label{EUB}
  (\tilde{\cL}_0x)(t)+(\tilde{\cL}_1a)(t)+(\tilde{\cL}_2b)(t)+(\tilde{\cL}_3\la_{a, b})(t)+(\tilde{\cL}_4\ga_{a, b})(t)=b(t).  
\end{eqnarray}

\begin{remark}
   Under an additional condition, we can prove that $\begin{pmatrix} \cL_3& \cL_4\\
   \tilde{\cL}_3& \tilde{\cL}_4\end{pmatrix}$  is invertible. In this case, we can uniquely solve $\la_{a, b}$ and $\ga_{a, b}$  from (\ref{EXB})  and (\ref{EUB})  in term of $a$ and $b$,  and then,  Problem (MF-SLQ)  will be converted to an  SLQ control problem with deterministic control variables $a$ and $b$. However, as we  have  seen from  Theorem \ref{BUAB},   the uniqueness of $\la_{a, b}$ and $\ga_{a, b}$ are not important.
\end{remark}

After we obtain the optimal $u^*_{a, b}$ for Problem 1, we finally consider the problem  with  the objective function 
   \begin{eqnarray}
  J_{a,b}( u^*_{a, b}) 
  &=&\EE\bigg\{\big<GX(T), X(T)\big>
+\int_0^T \Big(\big<Q(t)X(t), X(t)\big> +  \big<Q_1(t)a(t), a(t)\big>\nonumber\\
&&\ \ \ +\big<R(t)u^*_{a, b}(t), u^*_{a, b}(t)\big>+\big<R_1(t)b(t), b(t)\big> \Big)dt\bigg\}.
 \end{eqnarray}
where the  state process $X$ is the solution of equation (\ref{EXcsate}) with $u$ there replaced by $u^*_{a, b}$.

$\bf{Problem}$$\ \bm{2}$.  \label{oop}
  Find  $a^*\in L^2(0,T; \RR^{n})$ and  $b^*\in L^2(0,T; \RR^{m})$, such that 
  \begin{eqnarray*}
   J_{a^*, b^*}\big(u^*_{a^*, b^*}\big)
    =\inf_{a\in L^2(0,T; \RR^n), b\in L^2(0,T; \RR^m)}J_{a,b}( u^*_{a, b}) .
  \end{eqnarray*}
  
 Combining  the state feedback form of $u^*_{a, b}$ for Problem 1 given by  (\ref{OPu'''}) (with $\la$ and  $\ga$  replaced by $\la_{a, b}$ and $\ga_{a, b}$, respectively) and  the relationship between $ (a, b)$ and $(\la_{a, b}, \ga_{a, b})$  described through (\ref{EXB}) and (\ref{EUB}) implies that the unique optimal control $u^*_{a, b}$ for Problem 1 can be written as 
 \begin{eqnarray}\label{OPBU}
    u^*_{a, b}&=& (\cP_1x)+(\cK_1 a)+(\cN_1 b).
 \end{eqnarray}
 Further, the state process $X$ for Problem 2 can be represented as 
  \begin{eqnarray}\label{OPBX}
   X= (\cP_2x)+(\cK_2 a)+(\cN_2 b), \ and  \   X(T)= \cP_3x+\cK_3 a+\cN_3 b,  
  \end{eqnarray}
  where $\cP_i, \cK_i$ and $\cN_i$ are linear operators, $i=1,2,3.$
 \begin{theorem}\label{LEAB}
     Assume that (A1)-(A3) hold. Then $a^*$ and $ b^*$   are optimal for Problem 2 if and only if
  \begin{eqnarray}\label{Cab}
   \big( \cT_1^* R\cT_1+L+\cT_2^*Q\cT_2+\cT_3^* G\cT_3\big) \begin{pmatrix}
a^*\\
b^*
\end{pmatrix}+\big( \cT_1^* R\cP_1+\cT_2^*Q\cP_2+\cT_3^* G\cP_3\big) x=0,   
  \end{eqnarray}
where $L=\begin{pmatrix}Q_1 & \!0 &\\0  & \!R_1\end{pmatrix}$, $\cT_1=(\cK_1, \ \cN_1)$, $\cT_2=( \cK_2,\ \cN_2 )$, $\cT_3=(
      \cK_3,\ \cN_3 
   )$ and $\cT_i^*$ are adjoint operators of $\cT_i$ for $i=1,2,3.$
 \end{theorem}
 
Summing up the results above, we present the main theorem of this paper.
\begin{theorem}\label{FTH}
 Suppose that  (A1)-(A3) hold. Then, the optimal control $u^*$ of Problem (MF-SLQ)   as a feedback function   of  the optimal state $X^*$ is given by 
\[u^*=\Theta^{-1} N X^*+\Theta^{-1} M,\]
 where  $(X^*, M)$  is determined  by   (\ref {MPAO1}) and  (\ref{X0})  (with  $\la$ and  $\ga$ 
  replaced by any solution of equations (\ref{EXB}-\ref{EUB})  and       $ a$ and $b$ 
   replaced by  $a^*$ and $b^*$, respectively) while the deterministic process $a^*$   and $b^*$ are characterized by (\ref{Cab}).
\end{theorem}

\section{A preliminary stochastic maximum principle} \label{sec4}
\setcounter{equation}{0}
\renewcommand{\theequation}{\thesection.\arabic{equation}}

In this section, we will prove Theorems \ref{UST} and  \ref{pMSP} in detail by  using Mazur's theorem and by the convex variation principle, respectively.

{\it{The proof of Theorem \ref{UST}.}}
Similar to the proof of  Theorem 5.2 in  \cite{YZ1999}, using Mazur's theorem, we can obtain  the existence of an optimal control.  The uniqueness of the optimal control follows from the strict convexity of the cost functional $J(u)$.

 In fact,  the uniqueness of the solution to the state equation  implies that there exist some  operators $\cN$ and $ \hat{\cN} $  
   from $\RR^n$  to $L^2_{\FF}(0,T; \RR^n) $ and  $L^2_{\FF_{T}}(\Omega; \RR^n)$,  respectively,  $\cO$ and $\hat{\cO}$  from  $\cU_{ad}$   to $L^2_{\FF}(0,T; \RR^n) $ and  $L^2_{\FF_{T}}(\Omega; \RR^n)$,  respectively,  such that (see page 2817 in \cite{Y2013})
 \[X=\cN x+\cO u, \ \ X_T=\hat{\cN}x+\hat{\cO} u.\]
 For  the expectation  functional  $\EE:L_{\FF}^2(0,T; \HH)\to L^2(0,T; \HH)$ with $\HH=\RR^n, \RR^m$, we have 
\[\big<\bar{Q}_1 \EE[X], \EE[X]\big>_{L^2}=\big<(\EE^*\bar{Q}_1\EE)X, \ X\big>_{L_{\mathbb{F}}^2} \  and \ \big<\bar{R}_1\EE[u], \EE[u]\big>_{L^2}=\big<(\EE^*\bar{R}_1\EE)u, \ u\big>_{L_{\mathbb{F}}^2}.\]

 Therefore, the cost function $J$ can be written as a quadratic function of $u$  with the  form 
\begin{eqnarray*}
    J( u)&=&\EE\big[\big<GX_T, \ X_T\big>\big]+\big<QX, \ X\big>_{L_{\mathbb{F}}^2}+\big<\bar{Q}_1\EE[X], \ \EE[X]\big>_{L^2}+\big<R u, \ u\big>_{L^2_{\mathbb{F}}}+\big<\bar{R}_1 \EE[u], \ \EE[u]\big>_{L^2}\\
    &=&\Big<\big(\hat{\cQ}^*G\hat{\cO}+\cO^*Q\cO+\cO^*(\EE^*\bar{Q}_1\EE)\cO+R+\EE^*\bar{R}_1\EE\big)u, u\Big>_{L^2_{\mathbb{F}}}\\
    &&\ \ + 2\Big<\big(\hat{\cO}^*G\hat{\cN}+\cO^*Q\cN+\cO^*(\EE^*\bar{Q}_1\EE)\cN\big)x, u\Big>_{L^2_{\mathbb{F}}}\\
    &&\ \ +\Big<\big(\hat{\cN}^*G\hat{\cN}+\cN^*Q\cN+\cN^*(\EE^*\bar{Q}_1\EE)\cN\big)x, x\Big>_{L^2_{\mathbb{F}}}. 
\end{eqnarray*}
  Assumption (A3) implies
\[\hat{\cQ}^*G\hat{\cO}+\cO^*Q\cO+\cO^*\EE^*Q_1\EE\cO+R+\EE^*R_1\EE\geq R\geq \delta I_m.\]
Hence,  for any $\la\in (0,1)$ and $u_1,\ u_2\in\cU_{ad}$, we have 
 \[\la J( u_1)+(1-\la) J( u_2)>J( \la u_1+(1-\la )u_2).\]
 This yields the uniqueness of the optimal control, and hence, finishes the proof of Theorem  \ref{UST}.
\qed

Let $(u^*, X^*) $ be the optimal pair that satisfies equation $(\ref{sate})$.  Let $u^\ep=u^*+\ep v,\; \ep \in\RR,\; v\in\cU_{ad}$   and  
let $X^\ep$ denote the state trajectory corresponding  to the control $u^\ep$. Then, $X^\ep$  satisfies the following equation 
\begin{equation}\label{rsate}
\left\{\begin{array}{ccl}
dX^\ep(t)&=&\big(AX^\ep+A_1\bar{X}^\ep+Bu^\ep+B_1\bar{u}^\ep\big)dt +\big(CX^\ep+C_1\bar{X}^\ep+Du+D_1\bar{u}^\ep \big)dW(t)\\
&&\ \ +\int_{U_0}\big(\al(z)X^\ep+\al_1(z)\bar{X}^\ep+\be(z)u^\ep+\be_1(z)\bar{u}^\ep\big)\tilde{N}(dt,dz)\\
X^\ep(0)&=&x, \ t\in[0,T].
\end{array}\right.
\end{equation}

Denote  $X^1(t):=\ep^{-1}\big(X^\ep(t)-X^*(t)\big)$. By the linearity,   it is clear that $X^1(t)$  is  independent of $\ep$ and satisfies
\begin{equation}\label{rsateJ}
\left\{\begin{array}{ccl}
dX^1(t)&=&\big(A X^1 +A_1 \bar{X}^1+B v +B_1 \bar{v} \big)dt +\big(C X^1 +C_1 \bar{X}^1 +D v +D_1 \bar{v}  \big)dW(t)\\
&&\ \ +\int_{U_0}\big(\al(z)X^1+\al_1(z)\bar{X}^1+\be(z)v+\be_1(z)\bar{v}\big)\tilde{N}(dt,dz)\\
X^1(0)&=&0, \ t\in[0,T].
\end{array}\right.
\end{equation}

 Now we are ready to present
 
{\it{The proof of Theorem  \ref{pMSP}.}} Let $u^*\in \cU_{ad}$ and $(Y^*,Z^*, K^*)$   be the adapted solution  of the second equation in (\ref{FBSDE11}). By the definition, $u^*$ is optimal if and only if 
\begin{align}
  J(u^*+\ep v) -J(u^*) \geq 0, \  \forall \ep \in\RR, \ \forall v\in \cU_{ad}.
\end{align}

Substituting for $u^\ep=u^*+\ep v$  and $X^\ep= X^*+\ep X^1$ into $J(u^\ep)$, for any $\ep\in\RR, v\in\cU_{ad}$,  we have 
\begin{eqnarray*}
0&\leq&   J(u^\ep)  -J(u^*)\\
&=&\ep^2\EE\bigg\{\big<G X^1(T), X^1(T)\big>+\int_0^T\Big(\big<Q X^1, X^1\big>+\big<Q_1 \bar{X}^1, \bar{X}^1\big>+\big<R v, v\big>+\big<R_1 \bar{v}, \bar{v}\big>\Big)dt\bigg\}\\
&&+2\ep \EE\bigg\{\big<G X^*(T), X^1(T)\big>+\int_0^T\Big(\big<Q X^*, X^1\big>+\big<Q_1\bar{X}^*, \bar{X}^1\big> +\big<Ru^*, v\big>+\big<R_1\bar{u}^*, \bar{v}\big>\Big)dt\bigg\}\\
&=&\ep^2J(v)+2\ep \EE\bigg\{\big<G X^*(T), X^1(T)\big>+\int_0^T\Big(\big<QX^*, X^1\big>+\big<Q_1\bar{X}^*, \bar{X}^1\big>\\
&&\ \ \ \ \ +\big<Ru^*, v\big>+\big<R_1\bar{u}^*, \bar{v}\big>\Big)dt\bigg\}.
\end{eqnarray*}
By the positive definiteness of the weighting matrices, it is clear that $J(v)\geq 0,\ \ \forall v\in\cU_{ad}$. Hence, the above inequality is equivalent to 
\begin{equation}\label{eq1''}
 \EE\bigg\{\big<G X^*(T), X^1(T)\big>+\int_0^T\Big(\big<Q X^*, X^1\big>+\big<Q_1\bar{X}^*, \bar{X}^1\big>+\big<Ru^*, v\big>+\big<R_1\bar{u}^*, \bar{v}\big>\Big)dt\bigg\} = 0.
\end{equation}
Integrating by part, we obtain 
\begin{eqnarray}\label{YX1}
 &&\EE[\big<GX^*(T), X^1(T)\big>]\nonumber\\ 
 &= &\EE\big[\big<Y^*(T), X^1(T)\big>\big]-\EE\big[\big<Y^*(0), X^1(0)\big>\big]\nonumber\\
 &=& \EE\bigg\{\int_0^T-\<A^\top Y^*+C^\top Z^*+\int_{U_0}\al^\top(z)K^*(z)\nu(dz)+QX^*,\ X^1\>dt\nonumber\\
 &&+\int_0^T-\<\EE\big[A_1^\top Y^*+C_1^\top Z^*+\int_{U_0}\al_1^\top(z)K^*(z)\nu(dz)+Q_1\bar{X}^*\big],  \ X^1(t)\>dt\bigg\}\nonumber\\\
 &&+\EE\bigg\{\int_0^T\<Y^*, \ AX^1+A_1\bar{X}^1+Bv+B_1\bar{v} \>dt \bigg \}\nonumber\\
 &&+\EE\bigg\{\int_0^T\<Z^*,\ CX^1+C_1\bar{X}^1+Dv+D_1\bar{v}\>dt\bigg\}\nonumber\\
&&+\EE\bigg\{\int_0^T\int_{U_0}\big<K^*(z),\  \al(z)X^1+\al_1(z)\bar{X}^1+\be(z)v+\be_1\bar{v} \big>\nu(dz)dt\bigg\}\nonumber\\
&=&\EE\bigg\{\int_0^T\Big(\big<Y^*,\  Bv+B_1\bar{v}\big>+\big<Z^*,\  Dv+D_1\bar{v}\big>\Big)dt\nonumber\\
&&\ \ \  +\int_0^T\int_{U_0}\big<K^*(z), \ \be(z)v+\be_1(z)\bar{v}\big>\nu(dz)dt\bigg\}\nonumber\\
&&\ \ \ -\EE\bigg\{\int_0^T\Big(\big<QX^*, \ X^1\big>+\big<\EE[Q_1\bar{X}^*], \ X^1\big>\Big)dt\bigg\}
\end{eqnarray}
Plugging (\ref{YX1}) into (\ref{eq1''}) we have
\begin{eqnarray}\label{CD'}
0&=&\EE\bigg\{\int_0^T\<B^\top Y^*+D^\top Z^*+\int_{U_0}\be^\top(z)K^*(z)\nu(dz)+Ru^*, \ v\>dt\bigg\}  \nonumber\\
&& + \EE\bigg\{\int_0^T\<B_1^\top Y^*+D_1^\top Z^*+\int_{U_0}\be_1^\top(z)K^*(z)\nu(dz)+R_1\bar{u}^*, \ \bar{v}\>dt\bigg\}\nonumber\\
&=& \EE\bigg\{\int_0^T\<B^\top Y^*+D^\top Z^*+\int_{U_0}\be^\top(z)K^*(z)\nu(dz)+R u^*, \ v\>dt\bigg\}  \nonumber\\
&& + \EE\bigg\{\int_0^T\<\EE\big[B_1^\top Y^*+D_1^\top Z^*+\int_{U_0}\be_1^\top(z)K^*(z)\nu(dz)+R_1 \bar{u}^*\big], \ v(t)\>dt\bigg\}\nonumber\\
&&-\EE\bigg\{\int_0^T\<\EE\big[B_1^\top Y^*+D_1^\top Z^*+\int_{U_0}\be_1^\top(z)K^*(z)\nu(dz)+R_1 \bar{u}^* \big], \ v\>dt\bigg\}\nonumber\\
&& + \EE\bigg\{\int_0^T\<B_1^\top Y^*+D_1^\top Z^*+\int_{U_0}\be_1^\top(z)K^*(z)\nu(dz)+R_1 \bar{u}^*, \ \bar{v}\>dt\bigg\}\nonumber\\
&=&\EE\int_0^T\big<\big(\cH_u+\EE[\cH_{\bar{u}}]\big)(t,X^*,\bar{X}^*,  u^*, \bar{u}^*, Y^*, Z^*, K^*), \ v\big>dt, \ \forall v\in \cU_{ad}.
\end{eqnarray}
The last equality above is based on the fact that
\begin{eqnarray*}
 &&\EE\bigg\{\int_0^T\<\EE\big[B_1^\top Y^*+D_1^\top Z^*+\int_{U_0}\be_1^\top(z)K^*(z)\nu(dz)+R_1 \bar{u}^*\big], \ v\>dt\bigg\}\nonumber\\
&=&\int_0^T\<\EE\big[B_1^\top Y^* +D_1^\top Z^*+\int_{U_0}\be_1^\top(z)K^*(z)\nu(dz)+R_1\bar{u}^*\big], \ \bar{v}\>dt\nonumber\\
&=&\EE\bigg\{\int_0^T\<B_1^\top Y^*+D_1^\top Z^*+\int_{U_0}\be_1^\top(z)K^*(z)\nu(dz)+R_1 \bar{u}^*, \ \bar{v}\>dt\bigg\}.
\end{eqnarray*}
Therefore,  (\ref{CD'}) implies that Theorem \ref{pMSP}  holds.
\qed

\begin{remark}
Note that  the optimality system  (\ref{FBSDE11}-\ref{SCUB})  is a fully coupled FBSDE. It is difficult to find the optimal control by sloving   (\ref{FBSDE11}-\ref{SCUB}) directly  due to the terms  $\EE[A_1^\top Y^*+C_1^\top Z^*]$ and $\EE[R_1u^*]$  in the adjoint equation and  stationary condition, respectively. However,  using the ELM method in the next section, we will find the procedure to obtain the optimal control. 
\end{remark}

\section{ Extened LaGrange multiplier method} \label{sec5}
\setcounter{equation}{0}
\renewcommand{\theequation}{\thesection.\arabic{equation}}

 In this section, we consider the  Problem 1  for the study of the  non-homogeneous SLQ control  problem with random coefficients and L\'evy driving noises,  which extends the results of Meng \cite{Meng2014}  for the homogeneous case.  

\subsection{ Optimize the control  variable  for fixed ELMs} \label{sec51}
\setcounter{equation}{0}
\renewcommand{\theequation}{\thesection.\arabic{equation}}

 In this su'bsection, we consider   Problem (Sub1).  First, we present the proofs of Lemmas  \ref{P1U1}, \ref{UN2} and  \ref{OPSY2}. 
Finally, we will decouple the FBSDE (\ref{OMCadjion}), which leads to Theorem \ref{TH2'}.

{\it{The proof of Lemma  \ref{P1U1}.}}
 Similar  to the proof of Theorem {\ref{UST},  for $x\in \RR^{n},\ a\in L^2(0,T; \RR^n),\ b\in L^2(0,T; \RR^n)$ and $u\in \cU_{ad}$  fixed, the  uniqueness of SDE (\ref{EXcsate})  implies  process $X$ can be written as a linear combination of $x, u, a$   and $b$.  Therefore, the cost function $J_{a, b}$ can be written as a quadratic form of $u$ with  positive coefficients  (in fact, the  positive coefficients is greater than or equal to $R$)  when $x, a$ and $b$ are fixed.  
 
 Furthermore, $\EE[X^u]=a$ can be written as $\cN u=k$, where $\cN$ is a linear operator from $\cU_{ad}$ to $L^2([0,T], \RR^n)$. Let $\cR=\begin{pmatrix}
   \cN\\
   \EE
 \end{pmatrix}$ and  $d_1=\begin{pmatrix}
   k\\
   b
 \end{pmatrix} $. Then the   constraints  $\EE[X^u]=a$ and $\EE[u]=b$ can be replaced by $\cR u=d_1$. Note that $d_1\in Range(\cR)$,   there exists
  $v\in \cU_{ad}$ such that $\cR v=d_1$ and $\cR(u-v)=0$. Hence, we have $u-v\in\Ker(\cR)$, which  implies that there exists $d\in \Ker(\cR)$  such that $u=v+d.$ So, the control variables in Problem 1  can be replaced through $v+d$ and without constraints.    
 
 Combining the above analyses, we can obtain  the existence of a unique optimal control $u^*_{a, b}$ for the constrained problem. 
\qed

{\it{The proof of Lemma  \ref{UN2}.}}
 Similar to the proof of Lemma \ref{P1U1} above, we obtain  the existence and uniqueness of the  optimal control $u_{a, b, \la, \ga}$ when there is no constraint while $x, a, b, \la$ and $\ga$ are fixed. 
\qed

{\it{The proof of Lemma  \ref{OPSY2}.}}
Let $(u_{a,b, \la, \ga}, X_{a,b, \la, \ga})) $ be the optimal pair of Problem (Sub1). For any fixed but arbitrary $\ep\in \RR$,  let $u^\ep_{a,b, \la, \ga}=u_{a,b, \la, \ga}+\ep v,\; v\in\cU_{ad}$   and  let $X^\ep_{a,b, \la, \ga}$ denoted the state trajectory  to the control $u^\ep_{a,b, \la, \ga}$.
Let  $X^1_{a, b, \la, \ga}(t):=\ep^{-1}\big(X^\ep_{a, b, \la, \ga}(t)-X_{a, b, \la, \ga}(t)\big)$. It is clear that $X^1(t)$ satisfies the variation equation:
\begin{equation}\label{EXVE2}
\left\{\begin{array}{ccl}
dX^1_{a,b, \la, \ga}(t)&=&\big(AX^1_{a,b, \la, \ga}+Bv\big)dt+\big(CX^1_{a,b, \la, \ga}+Dv \big)dW(t)\\
&&\ \ +\int_{U_0}\big(\al(z)X_{a,b, \la, \ga}^1+\be(z)v\big)\tilde{N}(dt,dz)\\
X^1_{a,b, \la, \ga}(0)&=&0, \ t\in[0,T].
\end{array}\right.
\end{equation}
 Note that  $u_{a, b, \la, \ga}$ is an optimal control if and only if
\begin{eqnarray}\label{J2V} 
0&\leq&   J_{a,b,\la, \ga}(u^\ep_{a,b,\la, \ga})  -J_{a,b,\la, \ga}(u_{a,b,\la, \ga})\nonumber\\
&=&\EE\bigg\{\big<G X_{a,b,\la, \ga}^\ep(T), X_{a,b,\la, \ga}^\ep(T)\big>-\big<G X_{a,b,\la, \ga}(T), X_{a,b,\la, \ga}(T)\big>\nonumber\\
&&\ \ +\int_0^T\Big(\big<QX_{a,b,\la, \ga}^\ep, X_{a,b,\la, \ga}^\ep\big>-\big<QX_{a,b,\la, \ga}, X_{a,b,\la, \ga}\big>+\big<Ru_{a,b,\la, \ga}^\ep, u_{a,b,\la, \ga}^\ep\big>\nonumber\\
&&\ \ -\big<Ru_{a,b,\la, \ga}, u_{a,b,\la, \ga}\big>+2\big<\la, X_{a,b,\la, \ga}^\ep\big>-2\big<\la, X_{a,b,\la, \ga}\big> +2\big<\ga,  u_{a,b,\la, \ga}^\ep\big>-2\big<\ga,  u_{a,b,\la, \ga}\big>\Big)dt\bigg\}\nonumber\\
&=& \ep^2 \EE\Big\{\big<GX^1_{a, b, \la, \ga}, \ X^1_{a, b, \la, \ga}\big>+\int_0^T\Big(\big<Q X^1_{a, b, \la, \ga}, X^1_{a, b, \la, \ga}\big>+\big<R v, v\big>\Big)dt\Big\}\nonumber\\
&&+2\ep\EE\bigg\{\big<G X_{a,b,\la, \ga}(T), X_{a,b,\la, \ga}^1(T)\big> +\int_0^T\Big(\big<QX_{a,b,\la, \ga}, X_{a,b,\la, \ga}^1\big>+\big<Ru_{a,b,\la, \ga}, v\big>\nonumber\\
&&\ \ \ \ +\big<\la, X_{a,b,\la, \ga}^1\big>+\big<\ga,  v\big>\Big)dt\bigg\}.
\end{eqnarray}
It is clear from $(A3)$ that 
\[\EE\Big\{\big<GX^1_{a, b, \la, \ga}, \ X^1_{a, b, \la, \ga}\big>+\int_0^T\Big(\big<Q X^1_{a, b, \la, \ga}, X^1_{a, b, \la, \ga}\big>+\big<R v, v\big>\Big)dt\Big\}\geq 0.\]

For any fixed $v$ and $x$, the right-hand term of (\ref{J2V})  is a quadratic polynomial of the variable $\ep$. Hence,  (\ref{J2V}) holds if and only if
\begin{eqnarray}\label{EP2H}
&&\EE\int_0^T\Big(\<QX_{a,b,\la, \ga}, X_{a,b,\la, \ga}^1\>+\<Ru_{a,b,\la, \ga}, v\>+\<\la, X_{a,b,\la, \ga}^1\>+\<\ga,  v\>\Big)dt\nonumber\\
&+&\EE\<G X_{a,b,\la, \ga}(T), X_{a,b,\la, \ga}^1(T)\> =0.
\end{eqnarray}

  Applying It\^o's formula to $\big<Y_{a,b, \la, \ga}, X_{a,b, \la, \ga}^1\big>$,  we have 
\begin{eqnarray*} 
&&d\big<Y_{a,b, \la, \ga}(t), \ X^1_{a, b, \la, \ga}(t)\big>\nonumber\\
&=&- \<A^\top Y_{a,b, \la, \ga}+C^\top Z_{a,b, \la, \ga}+\int_{U_0}\al^\top(z)K_{a,b, \la, \ga}(z)\nu(dz)+Q X_{a, b, \la, \ga}+\la \big), X^1_{a, b, \la, \ga}\>dt\nonumber\\
&&\ \ \ +\<Z_{a,b, \la, \ga}, X^1_{a, b, \la, \ga}\>dW(t) +\int_{U_0}\< K_{a,b, \la, \ga}(z), \ X^1_{a, b, \la, \ga}\>\tilde{N}(dt,dz)dt\nonumber\\
&& + \<Y_{a,b, \la, \ga}, \  AX^1_{a,b, \la, \ga}+Bv \>dt+ \<Y_{a,b, \la, \ga}, CX^1_{a,b, \la, \ga}+Dv \>dW(t)\nonumber\\
&&\ \ \ + \int_{U_0} \<Y_{a,b, \la, \ga},\al(z)X_{a,b, \la, \ga}^1+\be(z)v\>\tilde{N}(dt,dz) \nonumber\\
&& +\big<Z_{a,b, \la, \ga}(t), \ CX^1_{a,b, \la, \ga}+Dv  \big>dt+\int_{U_0}\big<K_{a,b, \la, \ga}(z),  \ \al(z)X_{a,b, \la, \ga}^1+\be(z)v\big>\nu(dz)dt.\nonumber\\
&&+ \int_{U_0}\big<K_{a,b, \la, \ga}(z),  \ \al(z)X_{a,b, \la, \ga}^1+\be(z)v\big>\tilde{N}(dt,dz).
\end{eqnarray*}
Hence, 
\begin{eqnarray}\label{EXYX1} 
&&\EE\Big[\big<GX_{a, b, \la, \ga}, \  X^1_{a, b, a, \ga}\big>\Big]\\
&=&\EE\int_0^T\bigg\{\big<-\la-Q X_{a, b, \la, \ga},  \  X^1_{a, b, a, \ga}\big>\nonumber\\
&&+\Big<B^\top Y_{a, b, \la, \ga}+D^\top Z_{a, b, \la, \ga}+\int_{U_0}\al^\top(z) K_{a, b, \la, \ga}(z)\nu(dz), \ v\Big> \bigg\}dt. \nonumber
\end{eqnarray}
Substituting (\ref{EXYX1}) into (\ref{EP2H}), we have (\ref{J2V}) holds if and only if
\begin{eqnarray*} 
0=\EE\int_0^T\Big<Ru_{a, b, \la, \ga} +B^\top Y_{a, b, \la, \ga} + D^\top Z_{a, b, \la, \ga}+\int_{U_0}\al^\top(z) K_{a, b, \la, \ga}(z)\nu(dz)+\ga, \ v\Big>dt, \ \forall v\in \cU_{ad},
\end{eqnarray*} 
which implies that Lemma \ref{OPSY2} holds.
 \qed

Now, we decouple the  FBSDE (\ref{OMCadjion}).  To this end, we try the following form
\begin{equation}\label{CP}
Y_{a,b, \la, \ga}=PX_{a,b, \la, \ga}+\phi,
\end{equation}
where 
\[dP(t)=P_1dt+\Lambda dW(t)+\int_{U_0}\Gamma(z)\tilde{N}(dt,dz),  \ P(T)=G,\]
and 
\[d\phi(t)=\phi_1dt+\psi dW(t)+\int_{U_0}\theta(z)\tilde{N}(dt,dz), \  \phi(T)=0,\]
with $P_1$ and $\phi_1$ being determined later.

 Applying   It\^o's formula for general L\'evy-type stochastic integrals to  $Y_{a,b, \la, \ga}(t)$ in (\ref{CP}),  we get 
\begin{eqnarray}\label{FORY}
dY_{a,b, \la, \ga}(t)&=&d\big(P(t)X_{a,b, \la, \ga}(t)+\phi(t)\big)\nonumber\\
&=&\Big(P_1dt+\Lambda dW(t)+\int_{U_0}\Gamma(z)\tilde{N}(dt,dz)\Big)X_{a,b, \la, \ga}\nonumber\\
&&\ +P\Big(\big(AX_{a,b, \la, \ga}+A_1a+Bu_{a,b, \la, \ga}+B_1b\big)dt\nonumber\\
&&\ +\big(CX_{a,b, \la, \ga}+C_1a+Du_{a,b, \la, \ga}+D_1b \big)dW(t)\nonumber\\
&&\ +\int_{U_0}\big(\al(z)X_{a,b, \la, \ga}+\al_1(z)a+\be(z)u_{a,b, \la, \ga}+\be_1(z)b\big)\tilde{N}(dt,dz)\Big)\nonumber\\
&&+\Lambda\big(CX_{a,b, \la, \ga}+C_1a+Du_{a,b, \la, \ga}+D_1b \big)dt\nonumber\\
&&\ +\int_{U_0}\Gamma(z)\big(\al(z)X_{a,b, \la, \ga}+\al_1(z)a+\be(z)u_{a,b}+\be_1(z)b\big)N(dt,dz)\nonumber\\
&&\ \ +\phi_1dt+\psi dW(t)+\int_{U_0}\theta(z)\tilde{N}(dt,dz)\nonumber\\
&=&\Big(\phi_1+P_1X_{a,b, \la, \ga}+P\big(AX_{a,b, \la, \ga}+A_1a+Bu_{a,b, \la, \ga}+B_1b\big) \nonumber\\
&&\ \ +\Lambda\big(CX_{a,b, \la, \ga}+C_1a+Du_{a,b, \la, \ga}+D_1b \big)\nonumber\\
&&\ \ +\int_{U_0}\Gamma(z)\big(\al(z)X_{a,b, \la, \ga}+\al_1(z)a+\be(z)u_{a,b, \la, \ga}+\be_1(z)b\big)\nu(dz)\Big)dt\nonumber\\
&&+\Big(\psi+\Lambda X_{a,b, \la, \ga} +P\big(CX_{a,b, \la, \ga}+C_1a+Du_{a,b, \la, \ga}+D_1b \big)\Big)dW(t)\nonumber\\
&&+\int_{U_0}\Big(\theta(z)+\Gamma(z)X_{a,b, \la, \ga}+P\big(\al(z)X_{a,b, \la, \ga}+\al_1(z)a+\be(z)u_{a,b, \la, \ga}+\be_1(z)b\big)\nonumber\\
&&\ \ \ \ \ \  + \Gamma(z)\big(\al(z)X_{a,b, \la, \ga}+\al_1(z)a+\be(z)u_{a,b}+\be_1(z)b\big)\Big)\tilde{N}(dt,dz).\nonumber
\end{eqnarray}
 Comparing with the second equation in (\ref{OMCadjion}), we have  
\begin{eqnarray}\label{yzku}
 &&-\Big(A^\top Y_{a,b, \la, \ga}+C^\top Z_{a,b, \la, \ga}+\int_{U_0}\al^\top(z)K_{a,b, \la, \ga}(z)\nu(dz)+QX_{a,b, \la, \ga}+\la\Big)\nonumber\\
 &=&\phi_1+P_1X_{a,b, \la, \ga}+P\Big(AX_{a,b, \la, \ga}+A_1a+Bu_{a, b, \la, \ga}+B_1b\Big)\nonumber\\
&& +\Lambda\Big(CX_{a,b, \la, \ga}+C_1a+Du_{a, b, \la, \ga}+D_1b \Big)\nonumber\\
&&+\int_{U_0}\Gamma(z)\Big(\al(z)X_{a,b, \la, \ga}+\al_1(z)a+\be(z)u_{a, b, \la, \ga} +\be_1(z)b\Big)\nu(dz), 
\end{eqnarray}
\begin{eqnarray}\label{ZU}
  Z_{a,b, \la, \ga} & =& 
P\Big(CX_{a,b, \la, \ga}+C_1a+Du_{a, b, \la, \ga}+D_1b\Big)+\psi+\Lambda X_{a,b, \la, \ga}
\end{eqnarray}
and 
\begin{eqnarray}\label{KU}
  K_{a,b, \la, \ga}(z)
  &=&\big(P+\Gamma(z)\big)\al(z)X_{a,b, \la, \ga}+\big(P+\Gamma(z)\big)\be(z) u_{a, b, \la, \ga}+\big(P+\Gamma(z)\big)\al_1(z)a \nonumber\\
  &&\ \ \  +\big(P+\Gamma(z)\big)\be_1(z)b+\theta(z)+\Gamma(z)X_{a,b, \la, \ga}.
\end{eqnarray}
Now we calculate $u_{a, b, \la, \ga}$.  Plugging (\ref{CP}),  (\ref{ZU})  and (\ref{KU}) into (\ref{HST})
we obtain 
\begin{eqnarray*}
0&=&\Big(D^\top P  D +\int_{U_0}\be^\top(z) \big(P +\Gamma(z)\big)\be(z) \nu(dz)+R \Big)u_{a,b, \la, \ga} \nonumber\\
&&+\Big(B^\top  P +D^\top \Lambda +D^\top  P  C  +\int_{U_0}\big(\be^\top(z)\Gamma(z)+\be^\top(z)(P+\Gamma(z)) \al(z)\big)\nu(dz)\Big) X_{a,b, \la, \ga}\nonumber\\
&&+\Big(B^\top \phi+D^\top \psi+D^\top P C_1a+D^\top P D_1 b+\ga\nonumber\\
&&\ \ \ \ +\int_{U_0}\big(\be^\top \theta+\be^\top (P+\Gamma(z)) \al_1 a+\be^\top (P+\Gamma(z)) \be_1 b\big)\nu(dz)\Big)\nonumber\\
&&:=-\Theta u_{a,b, \la, \ga}+NX_{a,b, \la, \ga}+M.
\end{eqnarray*}
By  Lemma \ref{BSRE}  we know that  $\Theta$ is invertible,  the  optimal control can be represented  as 
\begin{eqnarray}\label{uad}
 u_{a,b, \la, \ga}=\Theta^{-1}(NX_{a,b, \la, \ga}+M).   
\end{eqnarray}
 where $\Theta, N$ are defined in Section 3 and 
\begin{eqnarray*}
M&=&B^\top \phi+D^\top \psi+D^\top P C_1a+D^\top P D_1 b+\ga\nonumber\\
&&\ \ \ \ +\int_{U_0}\Big(\be^\top \theta+\be^\top \big(P+ \Gamma(z)\big) \al_1 a+\be^\top \big(P+ \Gamma(z)\big) \be_1 b\Big)\nu(dz).
\end{eqnarray*}

To determine the equations for $P$ and $\phi$, we combine (\ref{yzku}),  (\ref{ZU}),  (\ref{KU}) and (\ref{uad}) to get 
\begin{eqnarray*}
 && \bigg\{P_1+A^\top P+P A+C^\top\Lambda+\Lambda C+\int_{U_0}\big(\al^\top(z)\Gamma(z)+\Gamma(z)\al(z)\big)\nu(dz)
 \nonumber\\
  &&\ \ +C^\top P  C
 +\int_{U_0}\big(\al^\top(z) P \al(z)+\al^\top(z) \Gamma(z) \al(z)\big)\nu(dz)+Q+N^\top\Theta^{-1}N\bigg\}X_{a,b, \la, \ga} \nonumber\\
  &&+\bigg\{ \phi_1+A^\top \phi+C^\top \psi+\int_{U_0}\al^\top(z)\theta(z)\nu(dz)+\la\nonumber\\
  &&\ \ \ +\Big( C^\top P C_1+P A_1+\Lambda C_1+\int_{U_0}\big(\al^\top(z) (P+\Gamma(z))\al_1(z)+\Gamma(z)\al_1(z)\big)\nu(dz)\Big)a\nonumber\\
  &&\ \ \ +\Big( C^\top P D_1+P B_1+\Lambda D_1+\int_{U_0}\big(\al^\top(z) (P+\Gamma(z))\be_1(z)+\Gamma(z)\be_1(z)\big)\nu(dz)\Big)b \nonumber\\
  &&\ \ \ +\Big( C^\top P D+P B+\Lambda D+\int_{U_0}\big(\al^\top(z) (P+\Gamma(z))\be(z)+\Gamma(z)\be(z)\big)\nu(dz)\Big)\Theta^{-1}M\bigg\}
  \\
  &&=0.
\end{eqnarray*}
Setting the coefficient in front of $X_{a,b, \la, \ga}$ to be $0$, we have
\begin{eqnarray*}
 P_1&=& -\bigg\{A^\top P+P A+C^\top\Lambda+\Lambda C+C^\top P  C
 +Q+N^\top\Theta^{-1}N\nonumber\\
  &&+\int_{U_0}\Big(\al^\top(z)\Gamma(z)+\Gamma(z)\al(z)+\al^\top(z)\big( P \Gamma(z)\big)\al(z)\Big)\nu(dz)   \bigg\} 
\end{eqnarray*}
 and 
 \begin{eqnarray*}
  \phi_1&= & -\bigg\{A^\top \phi+C^\top \psi+\int_{U_0}\al^\top(z)\theta(z)\nu(dz)+\la\nonumber\\
  &&\ \ \ +\Big( C^\top P C_1+P A_1+\Lambda C_1+\int_{U_0}\big(\al^\top(z) ( P+\Gamma(z) )\al_1(z)+\Gamma(z)\al_1(z)\big)\nu(dz)\Big)a\nonumber\\
  &&\ \ \ +\Big( C^\top P D_1+P B_1+\Lambda D_1+\int_{U_0}\big(\al^\top(z) (P+\Gamma(z))\be_1(z)+\Gamma(z)\be_1(z)\big)\nu(dz)\Big)b \nonumber\\
  &&\ \ \ +\Big( C^\top P D+P B+\Lambda D+\int_{U_0}\big(\al^\top(z) (P+\Gamma(z))\be(z)+\Gamma(z)\be(z)\big)\nu(dz)\Big)\Theta^{-1}M\bigg\}.
 \end{eqnarray*}
Hence, $P(t)$  satisfies the SREJ (\ref{RP})  and $\phi(t)$ satisfies the linear BSDE (\ref{Lphi}) with jump.

With the solvability of the optimal system (\ref{OMCadjion})  and the Riccati equation (\ref{RP}), which have been obtained in  Lemmas \ref{LOP1} and \ref{BSRE}, in the following we will show that the equation (\ref{Lphi}) is uniquely solvable.

{\it{The proof of Lemma   \ref{LBSDER}.}} 
Let $(X_{a, b, \la, \ga},  Y_{a, b, \la, \ga}, Z_{a, b, \la, \ga}, K_{a, b, \la, \ga}(z))$ be the unique solution of (\ref{OMCadjion}). Define 
\[ \phi=Y_{a, b, \la, \ga}-PX_{a, b, \la, \ga}, \]
\begin{eqnarray*}
\psi&=&Z_{a, b, \la, \ga}-\Lambda X_{a, b, \la, \ga}-P\(CX_{a,b, \la, \ga}+C_1a+D_1b\)\\
&&+PDR^{-1}\(B^\top Y_{a,b, \la, \ga}+D^\top Z_{a,b, \la, \ga}+\int_{U_0}\be^\top K_{a,b, \la, \ga}(z)\nu(dz)+\ga\)   ,
\end{eqnarray*}
and 
\begin{eqnarray*}
\theta(z)&=& K_{a, b,\la, \ga}(z,t)-\Gamma(z)X_{a, b, \la, \ga}-\big(P+\Gamma(z)\big)\Big(\al(z)X_{a,b, \la, \ga}+\al_1(z)a 
\\
&&\ \ \ \ \ \ -\be(z)R^{-1}\big(B^\top Y_{a,b, \la, \ga}+D^\top Z_{a,b, \la, \ga}+\int_{U_0}\be^\top K_{a,b, \la, \ga}(z)\nu(dz)+\ga \big)+\be_1(z)b\Big).
\end{eqnarray*}
 Applying It\^o's formula to $\phi$, it is easy to check that $(\phi, \psi, \theta(z))$  is a solution of (\ref{Lphi}). Also, we can check that  $(\phi, \psi, \theta(z)) \in L_{\FF}^2(0,T; \RR^n)\times L_{\FF}^{2, p}(0,T; \RR^n)\times \cL_{\FF}^{2, p}(0,T; \RR^n)$, $\forall p\in(1/2, 1)$  which follow from
 \begin{eqnarray*}
&& \EE\Big(\int_0^T|\Lambda(t) X_{a, b, \la, \ga}(t)|^2dt\Big)^{p}\\
&\leq& \EE\Big(\big(\sup_{t\leq T}|X_{a, b, \la, \ga}(t)|^{2}\big)^p\big(\int_0^T|\Lambda(t)|^2dt\big)^{p}\Big)\\
&\leq&  \Big(\EE\big(\sup_{t\leq T}|X_{a, b, \la, \ga}(t)|^2\big)\Big)^{p} \Big(\EE\big(\int_0^T|\Lambda(t)|^2dt\big)^{\frac{p}{1-p}}\Big)^{1-p}<\infty,
 \end{eqnarray*}
  and 
 \begin{eqnarray*}
&&\EE\Big(\int_0^T\int_{U_0}|\Gamma(z) X_{a, b, \la, \ga}(t)|^2\nu(dz)dt\Big)^{p}\\
&=& \EE\Big(\int_0^T|X_{a, b, \la, \ga}(t)|^2\int_{U_0}|\Gamma(z) |^{2}\nu(dz)dt\Big)^{p}\\
&\leq& \EE\Big(\big(\sup_{t\leq T}|X_{a, b, \la, \ga}(t)|^{2}\big)^p\big(\int_0^T\int_{U_0}|\Gamma(z) |^2\nu(dz)dt\big)^{p}\Big)\\
&\leq&  \Big(\EE\big(\sup_{t\leq T}|X_{a, b, \la, \ga}(t)|^{2}\big)\Big)^{p} \Big(\EE\big(\int_0^T\int_{U_0}|\Lambda(t)|^2\nu(dz)dt\big)^{\frac{p}{1-p}}\Big)^{{1-p}}<\infty.
 \end{eqnarray*}
 
Next we prove the uniqueness of the solution to  equation (\ref{Lphi}). For notational convenience, let  $X_{a, b, \la, \ga}, Y_{a, b, \la, \ga}, Z_{a, b, \la, \ga}$ be denoted as  $X, Y, Z$.   We take any solution $(\phi', \psi', \theta'(z))$  and define $X$ to be a solution to the first euqation of FBSDE (\ref{OMCadjion}) with  $Y=PX+\phi'$ and $Z, K(z)$ satisfy
\begin{eqnarray*} 
Z+PDR^{-1}D^\top Z&=&(PC+\Lambda)X-PDR^{-1}B^\top  Y-\int_{U_0}PDR^{-1}\be(z)^\top K(z)\nu(dz)-PR^{-1}\ga\\
&&+PC_1a +PD_1 b+ \psi',
\end{eqnarray*}
and
\begin{eqnarray*}
&&K(z)+(P+\Gamma(z))\be(z)R^{-1}\int_{U_0}\be^\top (z)K(z)\nu(dz)\\
&&=\big((P+\Gamma)\al(z)+\Gamma\big)X-(P+\Gamma) \be R^{-1}B^\top Y-(P+\Gamma)\be R^{-1}B^\top Z \\
&&-(P+\Gamma)R^{-1}\ga+(P+\Gamma)\al_1a +(P+\Gamma)\be_1 b+ \theta'(z).
\end{eqnarray*}
Then, for any $p\in(1/2,1)$, $(X, Y, Z, K)\in (S_{\FF}^2(0,T; \RR^n))^2\times L_{\FF}^{2, p}(0,T; \RR^n)\times \cL_{\FF}^{2, p}(0,T; \RR^n)$ is a solution to FBSDE (\ref{OMCadjion}). 
Next we prove that $Z\in L_\FF^2(0,T; \RR^n)$ and $K\in \cL_{\FF}^2(0,T; \RR^n)$.  For this purpose, 
we define the following stopping time
 \[\tau_k=\inf\bigg\{t>0; \int_0^t \Big(|Z|^2+\int_{U_0} K( z)|^2\nu(dz)\Big)ds\geq k\bigg\}.  \]

 Applying It\^o's formula to $ |Y(s)|^2$, we obtain that 
 \begin{eqnarray*}
\EE\big(|Y(T\wedge\tau_k)|^2-|Y(0)|^2\big)&=&-2\EE\int_0^{T\wedge\tau_k}\big<Y, \ A^\top Y+C^\top Z+\int_{U_0}\al^\top(z)K(z)\nu(dz)+QX+\la\big>ds\\
&&+\EE\int_0^{T\wedge\tau_k}\Big(Z^2+\int_{U_0}|K(z)|^2\nu(dz)\Big)ds.
\end{eqnarray*}
Hence 
\begin{eqnarray*}
&&\EE\int_0^{T\wedge\tau_k}\Big(Z^2+\int_{U_0}|K(z)|^2\nu(dz)\Big)ds\\
&\leq &\EE\big(|Y(T\wedge\tau_k)|^2-|Y(0)|^2\big)+K\EE\int_0^{T\wedge\tau_k}\big(|Y|^2+|X|^2+|\la|^2\big)ds
\\
&&+\frac{1}{2}\EE\int_0^{T\wedge\tau_k}|Z|^2ds+\frac{1}{2}\EE\int_0^{T\wedge\tau_k}\int_{U_0}|K(z)|^2\nu(dz)ds
\end{eqnarray*}
 Taking $k\to \infty$, we have 
 \begin{eqnarray*}
\EE\int_0^{T}\Big(Z^2+\int_{U_0}|K(z)|^2\nu(dz)\Big)ds&\leq &\EE\big(|Y(T)|^2-|Y(0)|^2\big)+K\EE\int_0^{T}\big(|Y|^2+|X|^2+|\la|^2\big)ds\\
&<&\infty.
\end{eqnarray*}
Therefore,  for any $t\in[0,T]$,  $(X, Y, Z, K)\in (S_{\FF}^2(0,T; \RR^n))^2\times L_{\FF}^{2}(0,T; \RR^n)\times \cL_{\FF}^{2}([0,T]\times U_0; \RR^n)$ is the unique  solution to FBSDE (\ref{OMCadjion}). Then uniqueness of $(\phi, \psi, \theta)$  follows from the uniqueness of the solution to (\ref{OMCadjion}). 
\qed

{\it{The proof of Theorem  \ref{TH2'}.}} 
Based on  Lemmas  \ref{BSRE} and \ref{LBSDER}, and (\ref{uad})
we can immediately derive the form of the optimal control $u_{a, b, \la, \ga}$  of Problem (Sub1). Substituting  (\ref{OPu'''})
into the state equation (\ref{EXcsate}), we obtain the corresponding optimal trajectory $X_{a, b, \la, \ga}$, which satisfies equation (\ref{X0}).

 \qed
 
\subsection{Optimize the ELMs} \label{sec52}
\setcounter{equation}{0}
\renewcommand{\theequation}{\thesection.\arabic{equation}}

 In this section, we consider Problem (Sub2) by presenting the proof of  Lemma {\ref{DMSP}}.

Let $(\la_{a,b}, \ \ga_{a,b})$ be the optimal control of Problem (Sub2). Let $$( X_{a,b, \la_{a, b}, \ga_{a, b}}, Y_{a,b, \la_{a, b}, \ga_{a, b}}, Z_{a,b, \la_{a, b}, \ga_{a, b}}, K_{a,b, \la_{a, b}, \ga_{a, b}}(z))$$
be the corresponding  optimal trajectory.  Let $( \al^1, \ \be^1)$ be a given process such that $( \la^\ep, \ \ga^\ep) \equiv ( \lambda_{a,b}+\ep \lambda^1, \ \gamma_{a,b}+\ep \gamma^1)\in L^2(0,T; \RR^{n+m})$, here $\ep \in[0, 1]$.  By optimality, we have
\begin{eqnarray}\label{JEP}
\tilde{J}_{a.b}\big(\lambda_{a,b}, \gamma_{a,b}\big)\geq \tilde{J}_{a,b}\big(\lambda^\ep , \gamma^\ep\big).
\end{eqnarray}
Let $\big(X_{a,b, \la, \ga}^\ep, Y_{a,b, \la, \ga}^\ep, Z_{a,b, \la, \ga}^\ep, K_{a,b, \la, \ga}^\ep(z)\big)$  be the state process with control
 $(\lambda^\ep, \ \gamma^\ep )$.

 Consider   the following variation FBSDE
  \begin{equation}\label{varia}
\left\{\begin{array}{ccl}
 dX^1(t)&=& \Big( AX^1-BR^{-1}\big( B^\top Y^1+D^\top Z^1+\int_{U_0}\be^{\top}(z)K^1(z)\nu(dz)+\gamma^1\big)\Big)dt\\
 &&+\Big(CX^1-DR^{-1}\big(B^\top Y^1+D^\top Z^1 +\int_{U_0}\be^{\top}(z)K^1(z)\nu(dz)+\gamma^1\big)\Big)dW(t)\\
 &&+\int_{U_0}\Big(\al(z)X^1-\be(z)R^{-1}\big(B^\top Y^1+D^\top Z^1+\int_{U_0}\be^{\top}(z)K^1(z)\nu(dz)\\
&&\ \ \ +\gamma^1\big)\Big)\tilde{N}(dt,dz)\\
dY^1(t)&=&-\Big(A^\top Y^1+C^\top Z^1+\int_{U_0}\al^\top(z)K^1(z)\nu(dz)+Q X^1+\lambda^1(t)\Big)dt\\
&&\ \ +Z^1(t)dW(t)+\int_{U_0}K^1(t,z)\tilde{N}(dt,dz)\\
X^1(0)&=&0,\  Y^1(T)=GX^1(T),  \ t\in[0,T],
\end{array}\right.
\end{equation}

  \begin{lem}\label{6.3}
Let (A1)-(A3)  hold. Then for any $t\in[0,T]$, 
\begin{eqnarray*}
&& X_{a,b,\la, \ga}^\ep= \ep  X^1+X_{a,b, \la_{a, b}, \ga_{a, b}}; \ Y_{a,b,\la, \ga}^\ep= \ep  Y^1+Y_{a,b, \la_{a, b}, \ga_{a, b}};\\
&& Z_{a,b,\la, \ga}^\ep= \ep  Z^1+Z_{a,b, \la_{a, b}, \ga_{a, b}}; \ K_{a,b, \la, \ga}^\ep= \ep   K^1+ K_{a,b, \la_{a, b}, \ga_{a, b}}. 
\end{eqnarray*}
  \end{lem}
  
  Now we are ready to present
  
{\it{The proof of Lemma  \ref{DMSP}:}}
If $( \lambda_{a, b}, \gamma_{a, b})$ is optimal, then  (\ref{JEP})  hold. Due to the expansion of cost functional, we obtain
\begin{eqnarray}\label{DOJEP-J}
 0&\geq &\lim_{\ep \to 0}\ep^{-1}\big(\tilde{J}_{a,b}(\lambda^\ep, \gamma^\ep)-\tilde{J}_{a,b}(  \lambda_{a, b}, \gamma_{a, b})\big) \nonumber\\
 &\to&2\EE[\big<G X_{a,b, \la_{a, b}, \ga_{a, b}}(T), X^1(T)\big>] +2\EE\int_0^T\bigg\{\big<Q X_{a,b, \la_{a, b}, \ga_{a, b}}, \ X^1\big>+\big<X_{a,b, \la_{a, b}, \ga_{a, b}}-a, \ \la^1\big>
\nonumber\\
&&\ \ +\big< \la_{a, b}, \  \ X^1\big> +\big<B R^{-1}B^{\top}Y_{a,b, \la_{a, b}, \ga_{a, b}}, \ Y^1\big>+\big<D R^{-1} D^\top  Z_{a,b, \la_{a, b}, \ga_{a, b}}, \ Z^1 \big>
\nonumber\\
 &&\ \ +\big<R^{-1}\int_{U_0}\be^\top(z)K_{a,b, \la_{a, b}, \ga_{a, b}}(z)\nu(dz), \ \int_{U_0}\be^\top(z)K^1(z)\nu(dz) \big>\nonumber\\
 &&\ \ +\big<D R^{-1}B^\top Y^1, \  Z_{a,b, \la_{a, b}, \ga_{a, b}}\big>+\big<D R^{-1}B^\top Y_{a,b, \la_{a, b}, \ga_{a, b}}, \  Z^1\big>\nonumber\\
 &&\ \ +
 \big<R^{-1}B^\top Y^1, \ \int_{U_0}\be^\top(z)K_{a,b, \la_{a, b}, \ga_{a, b}}(z)\nu(dz)\big>\nonumber\\
&& +\big<R^{-1}B^\top Y_{a,b, \la_{a, b}, \ga_{a, b}}, \ \int_{U_0}\be^\top(z)K^1(z)\nu(dz)\big>\nonumber\\
 &&\ \ +\<R^{-1}D^\top Z^1, \ \int_{U_0}\be^\top(z)K_{a,b, \la_{a, b}, \ga_{a, b}}(z)\nu(dz)\>\nonumber\\
 &&+\big<R^{-1}D^\top Z_{a,b, \la_{a, b}, \ga_{a, b}}, \ \int_{U_0}\be^\top(z)K^1(z)\nu(dz)\big>\nonumber\\
 &&\ \ -\big<b, \ \gamma^1\big>-\big<R^{-1}\ga_{a, b}, \ga^{1}\big>\bigg\}dt
\end{eqnarray}
 Applying It\^o's formula to $\big<Y_{a,b, \la_{a, b}, \ga_{a, b}}, \ X^1\big>$, we have 
\begin{eqnarray}\label{6.5}
&&d\big<Y_{a,b, \la_{a, b}, \ga_{a, b}}, X^1\big>\nonumber\\
&&=\big<dY_{a,b, \la_{a, b}, \ga_{a, b}}, \ X^1\big>+\big<Y_{a,b, \la_{a, b}, \ga_{a, b}}, \ dX^1\big>+d\big[Y_{a,b, \la_{a, b}, \ga_{a, b}}, X^1\big]_t\nonumber\\
&&=-\Big<A^\top Y_{a,b, \la_{a, b}, \ga_{a, b}}+C^\top Z_{a,b, \la_{a, b}, \ga_{a, b}}+Q X_{a,b, \la_{a, b}, \ga_{a, b}}+\lambda_{a, b}, \  X^1\Big>dt\nonumber\\
&&-\int_{U_0}\<\al^\top(z)K_{a,b, \la_{a, b}, \ga_{a, b}}(z), \  X^1\>\nu(dz)dt\nonumber\\
&&\ \ +\<Y_{a,b, \la_{a, b}, \ga_{a, b}}, \ AX^1-BR^{-1}\big( B^\top Y^1+D^\top Z^1+\int_{U_0}\be^{\top}(z)K^1(z)\nu(dz)+\gamma^1\big) \>dt\nonumber\\
&&\ \ +\Big<CX^1-DR^{-1}\big(B^\top Y^1+D^\top Z^1 +\int_{U_0}\be^{\top}(z)K^1(z)\nu(dz)+\gamma^1\big),  \ Z_{a,b, \la_{a, b}, \ga_{a, b}}\Big>dt\nonumber\\
&&\ \ +\int_{U_0}\Big<\al(z)X^1(t-)-\be(z)R^{-1}\big(B^\top Y^1+D^\top Z^1 +\gamma^1\big),  \ K_{a,b, \la_{a, b}, \ga_{a, b}}(z)\Big>\nu(dz)dt\nonumber\\
&&\ \ +\int_{U_0}\Big<-\be(z)R^{-1}\int_{U_0}\be^{\top}(z_1)K^1(z_1)\nu(dz_1),  \ K_{a,b, \la_{a, b}, \ga_{a, b}}(z)\Big>\nu(dz)dt\nonumber\\
&&\ \ + (\cdot\cdot\cdot)dW(t)+\int_{U_0}\cdot\cdot\cdot)\tilde{N}(dt.dz).
\end{eqnarray}
Hence, 
\begin{eqnarray*}
&&\EE[\big<GX_{a,b, \la_{a, b}, \ga_{a, b}}(T ), X^1(T)]\\
&=&\EE\int_0^Td\big<Y_{a,b, \la_{a, b}, \ga_{a, b}}(t), \ X^1(t)\big>\\
&=&\int_0^T\EE\bigg\{\big<-QX_{a, b, \la_{a, b}, \ga_{a, b}}-\la_{a, b}, \ X^1\big>+\big<Y_{a,b, \la_{a, b}, \ga_{a, b}}, \ -BR^{-1}\gamma^1\big>\ +\big<Z_{a,b, \la_{a, b}, \ga_{a, b}}, \ -DR^{-1}\gamma^1\big> \nonumber\\
&&\ \ +\int_{U_0}\big<K_{a,b, \la_{a, b}, \ga_{a, b}}(z), \ -\be(z)R^{-1}\gamma^1\big>\nu(dz)\\
&&\ \ - BR^{-1}\big<B^\top Y_{a,b, \la_{a, b}, \ga_{a, b}}+D^\top Z_{a,b, \la_{a, b}, \ga_{a, b}}+\int_{U_0}\be^\top(z)K_{a,b, \la_{a, b}, \ga_{a, b}}(z)\nu(dz), \ Y^1\big>\\
&&- DR^{-1}\big<B^{\top}Y_{a,b, \la_{a, b}, \ga_{a, b}}+D^\top Z_{a,b, \la_{a, b}, \ga_{a, b}}+\int_{U_0}\be^\top K_{a,b, \la_{a, b}, \ga_{a, b}}(z)\nu(dz), \ Z^1\big>\\
&&-\int_{U_0}\be(z) R^{-1}\Big<B^{\top}Y_{a,b, \la_{a, b}, \ga_{a, b}}+D^\top Z_{a,b, \la_{a, b}, \ga_{a, b}} +\int_{U_0}\be^\top(z) K_{a,b, \la_{a, b}, \ga_{a, b}}(z)\nu(dz) ,  \ K^1(z)\Big>\nu(dz)\bigg\}dt.
\end{eqnarray*}
Continuing with (\ref{DOJEP-J}), we have
\begin{eqnarray*}
 0&\geq&\lim_{\ep \to 0}\ep^{-1}\big(\tilde{J}_{a,b}(\lambda^\ep, \gamma^\ep)-\tilde{J}_{a,b}(\lambda_{a, b}, \gamma_{a, b})\big) \nonumber\\
 &=&
 \EE\int_0^T\bigg\{-\Big< R^{-1}\big( B^\top  Y_{a,b, \la_{a, b}, \ga_{a, b}} + D^\top  Z_{a,b, \la_{a, b}, \ga_{a, b}} +\ga_{a, b}\big)+b, \ \gamma^1\Big>\\
&&\ \ \ -\int_{U_0} \<R^{-1}\be^\top K_{a,b, \la_{a, b}, \ga_{a, b}}(z),\ga^1\>\nu(dz) +\big<X_{a,b, \la_{a, b}, \ga_{a, b}}-a, \ \lambda^1\big>\bigg\}dt. 
\end{eqnarray*}
 Together with (\ref{HST}), Lemma \ref{DMSP} is then proved.
\qed

{\it{The proof of Theorem  \ref{BUAB}.}}
As $u_{a, b, \la, \ga}$ is the optimal control for Problem (Sub1), we have 
\begin{eqnarray*}
J_{a, b, \la, \ga}(u_{a, b, \la, \ga})
\leq J_{a, b, \la, \ga}( u^*_{a, b})=J_{a, b}( u^*_{a, b})
\end{eqnarray*}
the above equation holds based on  $\EE[X^{u^*_{a, b}}-a]=0$ and $\EE[u^*_{a, b}]-b=0.$
For  EMLs $\la$ and $\ga$ such that   $\EE[X_{a, b, \la, \ga}]=a$ and $\EE[u_{a, b, \la, \ga}]=b$, we obtain
\begin{eqnarray*}
   J_{a, b, \la, \ga}(u_{a, b, \la, \ga})&=&J_{a, b} (u_{a, b, \la, \ga})+2\big<\la,\ \EE[X_{a, b, \la, \ga}]-a\big>_{L^2}+2\big<\ga,\ \EE[u_{a, b, \la, \ga}]-b\big>_{L^2}\\
   &=&J_{a, b} (u_{a, b, \la, \ga}).
\end{eqnarray*}
Hence, 
$$J_{a, b} (u_{a, b, \la, \ga})\leq J_{a, b}( u^*_{a, b}),$$
which becomes an equality due to the optimality of $u^*_{a, b}$.  The uniqueness of $u^*_{a,b}$  implies that $u_{a, b, \la, \ga}=u^*_{a, b}.$ 

Combining Theorem \ref{TH2'} and Lemma \ref{DMSP}, it is clear that the optimal pair $(u^*_{a, b}, X^*_{a, b})$ of Problem 1  is given by (\ref{OPu'''}) and (\ref{X0}), respectively  (with  $\la$ and  $\ga$  replaced by any solution of equations (\ref{DOS1}) and (\ref{DOS2})).  
\qed

\section{Optimal mean-field} \label{sec6}
\setcounter{equation}{0}
\renewcommand{\theequation}{\thesection.\arabic{equation}}

Finally, we will determine $(a, b)$ in this Section.  In other words, we present the proofs of {{Theorem}} \ref{LEAB} and Theorem \ref{FTH}. 
   For reader's convenience, we now restate Problem 2.  Let  $v=(a,\ b )^\top\in L^2(0,T; \RR^{n+m})$.   
   Recall that $u_{a, b}$ and $X$ have the  form   in (\ref{OPBU}) and (\ref{OPBX}). They are rewritten as
   \[u^*_{a,b}=(\cP_1 x)+(\cT_1 v)\]
   and 
   \[X=(\cP_2 x)+(\cT_2 v), \ and \ X(T)=\cP_3 x+\cT_3 v \]
   where $\cT_1=(\cK_1,\ \cN_1)$, $\cT_2=(\cK_2,\ \cN_2)$ and $\cT_3=(\cK_3,\ \cN_3)$.
   The cost function $J_{a, b}(u^*_{a, b})$ can be written as
   \begin{eqnarray*}
  J_{a,b}(u^*_{a, b}) 
  &=&\EE\bigg\{\big<G(\cP_3 x+\cT_3 v), \cP_3 x+\cT_3 v\big>\ \nonumber\\
&&\ \ \ +\int_0^T \Big(\Big<Q(t)\big((\cP_2 x)(t)+(\cT_2 v)(t)\big), (\cP_2 x)(t)+(\cT_2 v)(t)\Big> +  \big<L(t)v(t), v(t)\big>\nonumber\\
&&\ \ \ +\Big<R(t)\big((\cP_1 x)(t)+(\cT_1 v)(t)\big), (\cP_1 x)(t)+(\cT_1 v)(t)\Big>\Big)dt\bigg\}\nonumber\\
 &:=&J(v)
 \end{eqnarray*}
where 
$ L:=\begin{pmatrix}Q_1 & \!0 &\\0  & \!R_1\end{pmatrix}\in L^2_{\FF}(0,T; S_{n+m}).$

{\it{Proof of  {Theorem}  \ref{LEAB}.}} For any $(   a,\ b )^\top=v\in L^2(0,T; \RR^{n+m})$, we have 
\begin{eqnarray*}
  &&J_{a,b}(u^*_{a, b}) 
  =J(v)\nonumber\\  
  &=&\Big<\big( \cT_1^* R\cT_1+L+\cT_2^*Q\cT_2+\cT_3^* G\cT_3\big) v, v\Big>_{L^2}\nonumber\\
  &&+ 2\Big<\big( \cT_1^* R\cP_1+\cT_2^*Q\cP_2+\cT_3^* G\cP_3\big) x, v\Big>_{L^2}\nonumber\\
  &&+\Big<\big( \cP_1^* R\cP_1+\cP_2^*Q\cP_2+\cP_3^* G\cP_3\big) x, x\Big>_{L^2}.
\end{eqnarray*}
 Therefore,  $v^*=(   a^*,\ b^* )^\top$  is an optimal point  such that 
\begin{eqnarray*}
J( v^*)&=&\inf_{v\in L^2({0,T}; \RR^{n+m})}J(v)\\
&=&\inf_{a\in L^2(0,T;\RR^{n}), b\in L^2(0,T;\RR^m}J_{a, b}( u^*_{a, b})\\
&=&J_{a^*, b^*}(u^*_{a^*, b^*})
\end{eqnarray*}
if and only if 
 \[\big( \cT_1^* R\cT_1+L+\cT_2^*Q\cT_2+\cT_3^* G\cT_3\big) v^*+\big( \cT_1^* R\cP_1+\cT_2^*Q\cP_2+\cT_3^* G\cP_3\big) x=0.\]
 We then obtain the  conclusions of the {{Theorem}} \ref{LEAB}. 
\qed
  
 Now we are ready to present
 
{\it{Proof of Theorem  \ref{FTH}.}}
 It is easy to show that 
\begin{eqnarray*}
   \inf_{u\in\cU_{ad}}J( u)=\inf_{(a, b)\in L^2(0,T; \RR^{n+m})}\inf_{ u\in\cU_{ad}}\big\{J( u): \EE[X]=a, \EE[u]=b\big\}. 
\end{eqnarray*}

This justifies the decomposition of the original problem into Problems 1 and 2. By the 
unique solvability  of  Problems 1 and 2, it is clear that $u^*=u^*_{a^*, b^*}$.
Combining Theorems  \ref{BUAB}  and \ref{LEAB}, we finish the proof  of Theorem \ref{FTH}.  
\qed

\end{document}